\documentclass[smallextended]{svjour3}       
\smartqed  
\usepackage{graphicx}
\usepackage{epstopdf}
\usepackage{amsfonts}
\usepackage{amsbsy}
\usepackage{amssymb}
\usepackage{amsmath}
\usepackage{indentfirst}
\usepackage{geometry}
\usepackage{mathrsfs}
\usepackage{setspace}
\usepackage{multirow}
\usepackage{float}
\journalname{}
\begin{document}

\title{Distributionally Robust Joint Chance-Constrained Programming with Wasserstein Metric}
\author{Yining Gu   \and Yanjun Wang* }
\institute{This research was supported by National Natural Science Foundation of China (Grant No. 11271243).\\
\rule[0.4\baselineskip]{3.8cm}{0.5pt}\\
Yining Gu\at
School of Mathematics, Shanghai University of Finance and Economics, Shanghai, China\\
             E-mail: Yining\_Gu@163.com
           \and
          Yanjun Wang,  Corresponding author  \at
             School of Mathematics, Shanghai University of Finance and Economics, Shanghai, China \\
               E-mail: wangyj@mail.shufe.edu.cn
}
\date{Received: date / Accepted: date}
\maketitle
\begin{abstract}
In this paper, we develop an exact reformulation and a deterministic approximation for distributionally robust joint chance-constrained programmings $(\text{DRCCPs})$ with a general class of convex uncertain constraints under data-driven Wasserstein ambiguity sets. It is known that robust chance constraints can be conservatively approximated by worst-case conditional value-at-risk (CVaR) constraints. It is shown that the proposed worst-case CVaR approximation model can be reformulated as an optimization problem involving biconvex constraints for joint DRCCP. We then derive a convex relaxation of this approximation model by constructing new decision variables which allows us to eliminate biconvex terms. Specifically, when the constraint function is affine in both the decision variable and the uncertainty, then the resulting approximation model is equivalent to a tractable mixed-integer convex reformulation for joint binary DRCCP. Numerical results illustrate the computational effectiveness and superiority of the proposed formulations.
\end{abstract}

\setlength{\baselineskip}{2.5em}
\section{Introduction}
\subsection{Problem Setting}
In this paper,we study distributionally robust chance-constrained programmings $(\text{DRCCPs})$ of the form:
\begin{alignat}{2}
\min_{\boldsymbol{x}} \quad &\boldsymbol{c}^{\top}\boldsymbol{x} &\tag{1a}\\
\mbox{s.t.}\quad
&\boldsymbol{x}\in\mathit{S}, &\tag{1b}\\
&\!\!\inf_{\mathbb{P}\in\mathcal{P}}\mathbb{P}\left\{\boldsymbol{\xi}:f_{t}\left(\boldsymbol{x},\boldsymbol{\xi}\right)\geq 0, \forall t\in[T]\right\}\geq 1-\epsilon, &\tag{1c}
\end{alignat}
where $\boldsymbol{x} \in \mathbb{R}^n$ is a decision vector; the vector $\boldsymbol{c} \in \mathbb{R}^n$ represents the objective function
coefficients; the set $\mathit{S} \subseteq \mathbb{R}^n$ represents deterministic constraints on $\boldsymbol{x}$; $\boldsymbol{\xi}\in \mathbb{R}^m$ represents a $\text{m}$-dimensional random vector supported on $\mathrm{\Xi}\subseteq\mathbb{R}^m$; $\mathcal{P}$ is termed as an ``ambiguity set'' comprising all distributions that are compatible with the decision maker's prior information; the mapping $f_{t}\left(\boldsymbol{x},\boldsymbol{\xi}\right):\mathbb{R}^n \times \mathrm{\Xi}\rightarrow\mathbb{R}$ for any $t\in\left[T\right]:=\left\{1,2,\cdots,T\right\}$ represents a set of the uncertain constraints on $\boldsymbol{x}$. In addition, the distributionally robust chance constraint $(\text{DRCC})$ $(1\text{c})$ requires all $T$ uncertain constraints to be jointly satisfied for all the probability distributions from the ambiguity set with a probability of at least $1-\epsilon$, where $\epsilon\in\left(0,1\right)$ represents the risk level specified by the decision makers, and $\epsilon$ is often chosen to be small, e.g., 0.10 or 0.05. The problem $(1)$ is called a single or joint DRCCP if $T=1$ or $T>1$, respectively.\\
\indent We denote the feasible region induced by $(1\text{c})$ as
\begin{alignat}{2}
\mathit{Z_{D}}:&=\left\{\boldsymbol{x}\in \mathbb{R}^n:\inf_{\mathbb{P}\in\mathcal{P}}\mathbb{P}\left\{\boldsymbol{\xi}:f_{t}\left(\boldsymbol{x},\boldsymbol{\xi}\right)\geq 0, \forall t\in[T]\right\}\geq 1-\epsilon\right\}& \tag{2}\\
&=\left\{\boldsymbol{x}\in \mathbb{R}^n:\sup_{\mathbb{P}\in\mathcal{P}}\mathbb{P}\left\{\boldsymbol{\xi}:f_{t}\left(\boldsymbol{x},\boldsymbol{\xi}\right)< 0, \exists t\in[T]\right\}\leq \epsilon\right\}&\tag{3}\\
&=\left\{\boldsymbol{x}\in \mathbb{R}^n:\sup_{\mathbb{P}\in\mathcal{P}}\mathbb{E}_{\mathbb{P}}\left[\mathbb{I}_{\left\{f_{t}\left(\boldsymbol{x},\boldsymbol{\xi}\right)< 0, \exists t\in\left[T\right]\right\}}\left(\boldsymbol{\xi}\right)\right]\leq \epsilon\right\}.&\tag{4}
\end{alignat}
\indent Note that all results in the remainder of this paper are predicated on the following assumptions.\\
\noindent$\left(\mathbf{A1}\right)$ Each function $f_{t}\left(\boldsymbol{x},\boldsymbol{\xi}\right)$ is convex continuous in $\boldsymbol{\xi}$ for any fixed $\boldsymbol{x}$, and is concave continuous in $\boldsymbol{x}$ for any fixed $\boldsymbol{\xi}$.\\
\noindent$\left(\mathbf{A2}\right)$ The random vector $\boldsymbol{\xi}$ is supported on a nonempty closed convex set $\mathrm{\Xi}\subseteq\mathbb{R}^m$.\\
\indent In this paper, the ambiguity set employed in the distributionally robust formulation is the Wasserstein ball centered at the empirical distribution of the sample dataset. Moreover, we make the following assumption on the ambiguity set $\mathcal{P}$.\\
\noindent$\left(\mathbf{A3}\right)$ The Wasserstein ambiguity set $\mathcal{P}_{\mathit{W}}$ is defined as
\begin{alignat}{2}
\mathcal{P}_{\mathit{W}}=\left\{\mathbb{P}:\mathbb{P}\left\{\boldsymbol{\xi}\in\mathrm{\Xi}\right\}=1,
\mathit{W}\left(\mathbb{P},\mathbb{P}_{\boldsymbol{\tilde{\zeta}}}\right)\leq\delta\right\},\tag{5}
\end{alignat}
where the 1-Wasserstein metric is defined as
\begin{equation}\notag
\mathit{W}\left(\mathbb{P}_{1},\mathbb{P}_{2}\right)=\inf_{\mathbb{Q}}
\left\{ {\begin{array}{*{20}{c}}
{\int_{\mathrm{\Xi}\times\mathrm{\Xi}} \left \|\boldsymbol{\xi}_{1}-\boldsymbol{\xi}_{2}\right \|\mathbb{Q}\left(\mathit{d}\boldsymbol{\xi}_{1},\mathit{d}\boldsymbol{\xi}_{2}\right):}&\begin{array}{l}
\mathbb{Q}~\text{is a joint distributionally of}~\boldsymbol{\xi}_{1}~\text{and}~\boldsymbol{\xi}_{2}\\
\text{with marginals}~\mathbb{P}_{1}~\text{and}~\mathbb{P}_{2},\,\text{respectively}
\end{array}
\end{array}}\right\},
\end{equation}
for all distributions $\mathbb{P}_{1},\mathbb{P}_{2}\in\mathcal{M}(\mathrm{\Xi})$, where $\mathcal{M}(\mathrm{\Xi})$ contains all probability distributions $\mathbb{P}$ supported on $\mathrm{\Xi}$ with $\mathbb{E}_{\mathbb{P}} \left\|\boldsymbol{\xi}\right \|=\int_{\mathrm{\Xi}} \left \|\boldsymbol{\xi}\right \|\mathbb{P}\left(\mathit{d}\boldsymbol{\xi}\right)<\infty$.\\
\indent The 1-Wasserstein metric between $\mathbb{P}_{1}$ and $\mathbb{P}_{2}$, equipped with an arbitrary norm $\left \|\cdot\right \|$ on $\mathbb{R}^{m}$, represents the minimum transportation cost generated by moving the probability mass from $\mathbb{P}_{1}$ to $\mathbb{P}_{2}$. In $(5)$, $\mathbb{P}_{\boldsymbol{\tilde{\zeta}}}$ represents a discrete empirical distribution of $\boldsymbol{\tilde{\zeta}}$ with i.i.d. samples $\mathit{Z}=\left\{\boldsymbol{\zeta^{i}}\right\}_{i\in\left[N\right]}\subseteq\mathrm{\Xi}$ from the true distribution $\mathbb{P}^{\infty}$, i.e., its point mass function is $\mathbb{P}_{\boldsymbol{\tilde{\zeta}}}\left\{\boldsymbol{\tilde{\zeta}}=\boldsymbol{\zeta^{i}}\right\}=\frac{1}{N}$, and $\delta>0$ represents the Wasserstein radius.
\subsection{Literature Review}
\indent There are significant efforts on reformulations, approximations and convexity properties of DRCCP problems under various ambiguity sets. In particular, many approaches based on moments and statistical distances are commonly used to build ambiguity sets in DRCCP problems.\\
\indent We now review existing works on DRCCP problems with moment-based ambiguity sets [1-10]. It is well-known that more efforts have been made to derive tractable reformulations for single DRCCP with constraint functions that are affine in both the decision variable and the uncertainty. For instance, the authors in [1] demonstrated that with given first- and second-order moments, the set $\mathit{{Z}_{D}}$ for single DRCCP is equivalent to a tractable second-order conic representation. In [9], the authors developed tractable semidefinite programming for single DRCCP with given first- and second-order moments as well as the support of the uncertain parameters. In addition, the authors in [2] showed that the set $\mathit{{Z}_{D}}$ for single DRCCP is convex when $\mathcal{P}$ involves conic moment constraints or unimodality of $\mathbb{P}$. However, tractability results for joint DRCCP with affine uncertain constraints are very rare. It has been shown in [2] that the optimization problem over the set $\mathit{{Z}_{D}}$ is NP-hard in general. Thus, much of the earlier works derived deterministic approximations of the set $\mathit{{Z}_{D}}$ instead of developing its equivalent reformulations. For instance, in [7], with given first- and second-order moments, the authors proposed to mitigate the potential over-conservatism of the Bonferroni approximation by scaling each uncertain constraint with a positive number and then converted them into a single chance constraint. They also derived a conservative second-order conic programming approximation for any given scaler. In [9], with given first- and second-order moments, the authors derived a deterministic reformulation of the set $\mathit{{Z}_{D}}$, but it is not convex due to bilinear terms, which is naturally hard to solve. Besides, the authors in [11] provided several sufficient conditions under which the well-known Bonferroni approximation is exact and obtained its convex reformulation. On the other hand, there is very limited literature on DRCCP problems with a broader family of uncertain constraints. For instance, in [9], with given first- and second-order moments, the authors proved that single DRCCP amounts to a tractable semidefinite programming when the constraint function is either concave piecewise or (possibly non-concave) quadratic in the uncertainty. With any ambiguity set including convex moment constraints, the authors in [5] studied deterministic reformulations of the set $\mathit{{Z}_{D}}$ and its convexity properties when the constraint function is convex in the uncertainty and is concave in the decision variable. They investigated deterministic reformulations of such problems and proposed some conditions under which such deterministic reformulations are convex. In [6], the authors proved that with given first- and second-order moments, single DRCCP is equivalent to a robust optimization problem when the constraint function is quasi-convex in the uncertainty and is concave in the decision variable.\\
\indent Recently, there are many successful developments on DRCCP problems with Wasserstein ambiguity sets [12-17]. Based on [18, 19], it should be noted that the Wasserstein ambiguity sets offer powerful out-of-sample performance guarantees and asymptotic consistency, and also enable the decision makers to control conservativeness of the distribution uncertainty by tuning the radius. Most existing results on the tractability of DRCCP problems with Wasserstein ambiguity sets are restricted to the case of affine uncertain constraints. For instance, the authors in [12] derived exact mixed-integer conic reformulations for single DRCCP as well as joint DRCCP with right-hand side uncertainty. In [13], the author showed that joint DRCCP is mixed-integer representable by introducing big-M parameters and additional binary variables. He also derived tractable outer and inner approximations of the set $\mathit{{Z}_{D}}$. The authors in [14] provided exact reformulations for single DRCCP under discrete support and derived deterministic approximations for single DRCCP under continuous support. To the best of our knowledge, the case of a broader family of uncertain constraints is largely untouched. This kind of optimization problem is of great interest, because the uncertainty might not be inherently affine in many applications. One exception that we are aware of is [15], where the authors studied conditional value-at-risk (CVaR) approximation of single DRCCP with a general class of constraint functions. They considered many constraint functions that are convex in the decision variable and then replace the set $\mathit{{Z}_{D}}$ with convex CVaR approximation. They also presented tractable reformulation of the CVaR approximation when the constraint function is the maximum of functions that are affine in both the decision variable and the uncertainty, and the support of the uncertainty is a polyhedron. Moreover, when the constraint function is concave in the uncertainty, they showed that a central cutting surface algorithm [20, 21] for semi-infinite programmings can be used to compute an approximately optimal solution of the CVaR approximation of single DRCCP.
\subsection{Contributions, Structure, and Notations}
\indent In this paper, we consider joint DRCCP problems with a general class of convex
uncertain constraints under Wasserstein ambiguity sets. Specifically, our main contributions are summarized as below:\\
\indent $1$. We develop the worst-case CVaR approximation of the set $\mathit{Z}_{D}$ for more complicated and general joint DRCCP. It is proved that the resulting approximation is not convex in general since it involves biconvex constraints. We then derive a convex relaxation of the proposed biconvex approximation by constructing new decision variables which allows us to eliminate biconvex terms. It turns out that the proposed convex relaxation is essentially exact when there is a single uncertain constraint $\left(\text{i.e}., T=1\right)$.\\
\indent $2$. We transform the convex relaxation of the proposed biconvex approximation into tractable conic programming reformulation when the constraint function is quadratic convex in the uncertainty and concave in the decision variable, and the support of the uncertainty is polyhedron or ellipsoid.\\
\indent $3$. We demonstrate that for joint binary DRCCP, the proposed biconvex approximation admits a tractable mixed-integer convex reformulation when the constraint function is affine in both the decision variable and the uncertainty.\\
\indent The remainder of the paper is organized as follows. Section 2 presents an exact reformulation of the
set $\mathit{Z}_{D}$. Section 3 develops a biconvex approximation of the set $\mathit{Z}_{D}$ based on the worst-case CVaR constraints. Section 4 reports the numerical results to illustrate
the performance of the proposed models. Finally, Section 5 summarizes this paper.\\
\noindent $\textbf{Notation}$. The following notation is used throughout the paper. We use bold-letters $(\text{e.g}., \boldsymbol{x}, \boldsymbol{A})$ to denote vectors or matrices, and use corresponding non-bold letters to denote their components. We use $\boldsymbol{\xi}$ to denote a random vector and $\xi$ to denote a realization of $\boldsymbol{\xi}$. We let $\boldsymbol{e}_{n}$ be the all-one vectors with dimension $n$. Given a positive integer $n$, we let $\left[n\right]:=\left\{1,2,\cdots,n\right\}$, and $\mathbb{R}^{n}_{+}:=\left\{\boldsymbol{x}\in \mathbb{R}^{n}: x_{i}\geq 0, \forall i\in \left[n\right]\right\}$ and $\mathbb{R}^{n}_{++}:=\left\{\boldsymbol{x}\in \mathbb{R}^{n}: x_{i}>0, \forall i\in \left[n\right]\right\}$. We denote $\left(t\right)_{+}=\max\left\{t,0\right\}$ for any given real number $t$. We define the indicator function as $\mathbb{I}_{\mathit{A}}(\boldsymbol{\xi})=1$, if $\boldsymbol{\xi}\in \mathit{A}$; $=0$, otherwise. Similarlly, the characteristic function is defined as $\chi_{\mathit{A}}\left(\boldsymbol{\xi}\right)=0$, if $\boldsymbol{\xi}\in\mathit{A}$; $=\infty$, otherwise. Given a norm $\left \|\cdot\right \|$ on $\mathbb{R}^{n}$, the dual norm $\left \|\cdot\right \|_{*}$ is defined by $\left\|\boldsymbol{z}\right\|_{*}:=\sup_{\left\|\boldsymbol{\xi}\right\|\leq1}\boldsymbol{z}^{\top}\boldsymbol{\xi}$. The space of symmetric matrices of dimension $n$ is denoted by $\mathbb{S}^{n}$. For any two matrices $\boldsymbol{X}$, $\boldsymbol{Y}\in\mathbb{S}^{n}$, the relation $\boldsymbol{X}\succeq\boldsymbol{Y} (\boldsymbol{X}\succ\boldsymbol{Y})$ implies that $\boldsymbol{X}-\boldsymbol{Y}$ is positive semidefinite (positive definite). The space of positive semidefinite (or positive definite)
matrices of dimension $n$ is denoted by $\mathbb{S}^{n}_{+}~(\text{or}~\mathbb{S}^{n}_{{+}{+}})$. The inner product of two matrices $\boldsymbol{X},\boldsymbol{Y}\in\mathbb{R}^{m\times n}$ is denoted by
$\langle\boldsymbol{X},\boldsymbol{Y}\rangle=\text{tr}(\boldsymbol{X}\boldsymbol{Y})=
\sum_{i=1}^{m}\sum_{j=1}^{n}x_{ij}y_{ij}$. Additional notations will be introduced as needed.
\section{Exact Reformulation}
\indent In this section, we develop a deterministic reformulation of the set $\mathit{Z}_{D}$. The derivation of the exact reformulation utilizes strong duality result which is introduced in Lemma 1.\\
\indent We first review the well-known strong duality result from [18], which can be applied to formulate the worst-case chance constraint into its dual form, and indeed by the proof of strong duality theorem in [18], we present another equivalent dual reformulation in Lemma 1.
\begin{lemma}\label{1}
(Dual Reformulation, Theorem 4.2 in [18]) $\left(i\right)$ The uncertainty set $\mathrm{\Xi}\subseteq\mathbb{R}^m$ is convex and closed. $\left(ii\right)$  Let
$l\left(\boldsymbol{x},\boldsymbol{\xi}\right):=\max_{k\in\left[K\right]}l_{k}\left(\boldsymbol{x},\boldsymbol{\xi}\right)$ represent a pointwise maximum function, and the negative constituent functions $-l_{k},\forall k\in\left[K\right]$ are proper, convex and lower semi-continuous in $\boldsymbol{\xi}$ and assume that $l_{k},\forall k\in\left[K\right]$ is not identically $-\infty$ on
$\mathrm{\Xi}$. Then, for any given $\delta\geq 0$, the worst-case expectation $\sup_{\mathbb{P}\in\mathcal{P}_{\mathit{W}}}\mathbb{E}_{\mathbb{P}}\left[l(\boldsymbol{x},\boldsymbol{\xi})\right]$ equals the optimal value of the finite convex program
\begin{alignat}{2}
\inf_{\lambda,\boldsymbol{s},\boldsymbol{z}}~&\lambda\delta+\frac{1}{N}\sum_{i=1}^{N}{s_{i}}& \tag{6a}\\
\mbox{s.t.}\quad\!
&\sup_{\boldsymbol{\xi}\in\mathrm{\Xi}}\left[\boldsymbol{z}_{ik}^{\top}\boldsymbol{\xi}+l_{k}\left(\boldsymbol{x},\boldsymbol{\xi}\right)\right]\!-\boldsymbol{z}_{ik}^{\top}\boldsymbol{\zeta^{i}}\!\leq s_{i}, ~\forall i\in\left[N\right],\forall k\in\left[K\right], & \tag{6b}\\
&\left \|\boldsymbol{z}_{ik}\right \|_{*}\leq \lambda,~\forall i\in\left[N\right],\forall k\in\left[K\right],&\tag{6c}\\
&\lambda\!\geq 0,&\tag{6d}
\end{alignat}
where $\lambda$, $\boldsymbol{s}$ and $\boldsymbol{z}$ are decision variables, $\lambda\!\geq 0$ is the dual variable for Wasserstein metric constraint $\int_{\mathrm{\Xi}\times\mathrm{\Xi}} \left \|\boldsymbol{\xi}-\boldsymbol{\tilde{\zeta}}\right \|\mathbb{Q}\left(\mathit{d}\boldsymbol{\xi},\mathit{d}\boldsymbol{\tilde{\zeta}}\right)\leq\delta$, and $\left \|\boldsymbol{z}_{ik}\right \|_{*}$ is the dual norm.
\end{lemma}

\indent Next, we can represent the indicator function as a pointwise maximum of a finite number of concave functions by Lemma 2 due to [14, 18].
\begin{lemma}\label{2}
The indicator function $\mathbb{I}_{\left\{f_{t}\left(\boldsymbol{x},\boldsymbol{\xi}\right)< 0, \exists t\in\left[T\right]\right\}}\left(\boldsymbol{\xi}\right)$ can be rewritten as the pointwise maximum of a finite number of concave functions, which is defined as
\begin{equation}
\mathbb{I}_{\left\{f_{t}\left(\boldsymbol{x},\boldsymbol{\xi}\right)< 0, \exists t\in\left[T\right]\right\}}\left(\boldsymbol{\xi}\right)=\max\left\{1-\chi_{\left\{f_{1}\left(\boldsymbol{x},\boldsymbol{\xi}\right)< 0\right\}}\left(\boldsymbol{\xi}\right),\cdots,1-\chi_{\left\{f_{T}\left(\boldsymbol{x},\boldsymbol{\xi}\right)< 0\right\}}\left(\boldsymbol{\xi}\right),0\right\},\tag{7a}
\end{equation}
where for any $t\in\left[T\right]$
\begin{equation}\tag{7b}
\chi_{\left\{f_{t}\left(\boldsymbol{x},\boldsymbol{\xi}\right)< 0\right\}}\left(\boldsymbol{\xi}\right)=
\left \{
\begin{aligned}
&0,~~~\:\text{if}~f_{t}\left(\boldsymbol{x},\boldsymbol{\xi}\right)< 0\\
&\infty,~~\text{otherwise}
\end{aligned}
\right.
\end{equation}
which is the characteristic function of the open convex set defined by $f_{t}\left(\boldsymbol{x},\boldsymbol{\xi}\right)<0$.
\end{lemma}

\indent We develop an equivalent reformulation of the set $\mathit{Z}_{D}$ in the next theorem which applies to joint DRCCP with more general convex uncertain constraints.
\begin{theorem}\label{8}
The feasible set $\mathit{Z}_{D}$ is equivalent to
\begin{equation}\notag
~~~~~~~~~~~~~~~~~~\mathit{Z}_{D}=
\left\{{\begin{array}{*{20}{c}}{\boldsymbol{x} \in \mathbb{R}^n:}&\begin{array}{l}
\lambda\delta+\frac{1}{N}\sum_{i=1}^{N}{s_{i}}\leq\epsilon,\\
G_{f_{t}}\left(\boldsymbol{z}_{it},\eta_{it},\boldsymbol{x}\right)+1-\boldsymbol{z}_{it}^{\top}\boldsymbol{\zeta^{i}}-s_{i}\!\leq0,~\forall i\in\left[N\right],\forall t\in T\left(\boldsymbol{x}\right),\\
\left \|\boldsymbol{z}_{it}\right \|_{*}-\lambda\leq 0,~\forall i\in\left[N\right],\forall t\in T\left(\boldsymbol{x}\right),\\
\lambda\!\geq 0,~\eta_{it}\!\geq 0,s_{i}\geq0,~\forall i\in\left[N\right],\forall t\in T\left(\boldsymbol{x}\right),
\end{array}
\end{array}} \right\}\begin{array}{*{20}{c}}
{~~~~~~~~~~~~\left( 8a \right)}\\
{~~~~~~~~~~~~\left( 8b \right)}\\
{~~~~~~~~~~~~\left( 8c \right)}\\
{~~~~~~~~~~~~\left( 8d \right)}
\end{array}
\end{equation}
where $G_{f_{t}}\left(\boldsymbol{z}_{it},\eta_{it},\boldsymbol{x}\right)=
\sup_{\boldsymbol{\xi}\in\mathrm{\Xi}}\left[\boldsymbol{z}_{it}^{\top}\boldsymbol{\xi}
-\eta_{it}f_{t}\left(\boldsymbol{x},\boldsymbol{\xi}\right)\right]$ and $T\left(\boldsymbol{x}\right)=\left\{t\in\left[T\right]:\exists \boldsymbol{\xi}\in\mathrm{\Xi},\,f_{t}\left(\boldsymbol{x},\boldsymbol{\xi}\right)<0\right\}$.
\end{theorem}
{\it Proof}.~\;Note that
\begin{equation}\notag
\mathit{Z}_{D}:=\left\{\boldsymbol{x}\in \mathbb{R}^{n}:\sup_{\mathbb{P}\in\mathcal{P}}\mathbb{E}_{\mathbb{P}}\left[\mathbb{I}_{\left\{f_{t}\left(\boldsymbol{x},\boldsymbol{\xi}\right)< 0, \exists t\in\left[T\right]\right\}}\left(\boldsymbol{\xi}\right)\right]\leq \epsilon\right\}.
\end{equation}
\indent According to Lemma 2, the indicator function $\mathbb{I}_{\left\{f_{t}\left(\boldsymbol{x},\boldsymbol{\xi}\right)< 0, \exists t\in\left[T\right]\right\}}\left(\boldsymbol{\xi}\right)$ can be represented as the pointwise maximum as denoted in (7).\\
\indent Therefore, by Lemma 1, the left-hand side of the constraint defining $\mathit{Z}_{D}$ can be rewritten as
\begin{alignat}{2}
\inf_{\lambda,\boldsymbol{s},\boldsymbol{z}}~&\lambda\delta+\frac{1}{N}\sum_{i=1}^{N}{s_{i}}&\tag{9a}\\
\mbox{s.t.}\quad\!\!
&\sup_{\boldsymbol{\xi}\in\mathrm{\Xi}}\left[\boldsymbol{z}_{it}^{\top}\boldsymbol{\xi}+1-\chi_{\left\{f_{t}\left(\boldsymbol{x},\boldsymbol{\xi}\right)< 0\right\}}\left(\boldsymbol{\xi}\right)\right]-\boldsymbol{z}_{it}^{\top}\boldsymbol{\zeta^{i}}\!\leq s_{i},~\forall i\in\left[N\right],\forall t\in \left[T\right], &\tag{9b}\\
&\left \|\boldsymbol{z}_{it}\right \|_{*}\leq \lambda,~\forall i\in\left[N\right],\forall t\in \left[T\right],&\tag{9c}\\
&\lambda\!\geq 0,s_{i}\geq0,~\forall i\in\left[N\right].&\tag{9d}
\end{alignat}
\indent Then, the optimization problem $\sup_{\boldsymbol{\xi}\in\mathrm{\Xi}}\left[\boldsymbol{z}_{it}^{\top}\boldsymbol{\xi}+1-\chi_{\left\{f_{t}\left(\boldsymbol{x},\boldsymbol{\xi}\right)< 0\right\}}\left(\boldsymbol{\xi}\right)\right]$ in $(9\text{b})$ can be rewritten as
\begin{alignat}{2}
\sup_{\boldsymbol{\xi}\in\mathrm{\Xi}} \quad &\boldsymbol{z}_{it}^{\top}\boldsymbol{\xi}+1& \tag{10a}\\
\mbox{s.t.}\quad
&f_{t}\left(\boldsymbol{x},\boldsymbol{\xi}\right)<0,&\tag{10b}
\end{alignat}
for any $i\in\left[N\right]$ and $t\in T\left(\boldsymbol{x}\right)$.\\
\indent Hence, for any $i\in\left[N\right]$ and $t\in T\left(\boldsymbol{x}\right)$, we use Lagrangian duality result to reformulate the problem (10) as
\begin{alignat}{2}
\mathop{\sup}\limits_{\scriptstyle~~~\boldsymbol{\xi}\in\mathrm{\Xi}\hfill\atop
\scriptstyle~f_{t}\left(\boldsymbol{x},\boldsymbol{\xi}\right)<0\hfill}\!\!\!\boldsymbol{z}_{it}^{\top}\boldsymbol{\xi}+1=\!\!\!\!\!\!\!\mathop{\sup}\limits_{\scriptstyle~~~\boldsymbol{\xi}\in\mathrm{\Xi}\hfill\atop
\scriptstyle~f_{t}\left(\boldsymbol{x},\boldsymbol{\xi}\right)\leq0\hfill}\!\!\!\boldsymbol{z}_{it}^{\top}\boldsymbol{\xi}+1&=
\sup_{\boldsymbol{\xi}\in\mathrm{\Xi}}\inf_{\eta_{it}\geq 0}\left[\boldsymbol{z}_{it}^{\top}\boldsymbol{\xi}+1-\eta_{it}f_{t}\left(\boldsymbol{x},\boldsymbol{\xi}\right)\right]&\tag{11a}\\
&=\inf_{\eta_{it}\geq 0}\sup_{\boldsymbol{\xi}\in\mathrm{\Xi}}\left[\boldsymbol{z}_{it}^{\top}\boldsymbol{\xi}+1-\eta_{it}f_{t}\left(\boldsymbol{x},\boldsymbol{\xi}\right)\right]
,&\tag{11b}
\end{alignat}
where the first equality follows from the fact that for any given $\boldsymbol{x}$, $f_{t}\left(\boldsymbol{x},\boldsymbol{\xi}\right)$ and the objective function $\boldsymbol{z}_{it}^{\top}\boldsymbol{\xi}+1$ are both continuous in $\boldsymbol{\xi}$, $\mathrm{\Xi}$ is a nonempty closed set, so that we can replace $``<"$ by $``\leq"$ without effect on the supremum.\\
\indent Thus, for any $i\in\left[N\right]$ and $t\in T\left(\boldsymbol{x}\right)$, $(9\text{b})$ can be rewritten as
\begin{equation}\notag
\exists\eta_{it}\geq0,~\sup_{\boldsymbol{\xi}\in\mathrm{\Xi}}\left[\boldsymbol{z}_{it}^{\top}\boldsymbol{\xi}+1
-\eta_{it}f_{t}\left(\boldsymbol{x},\boldsymbol{\xi}\right)\right]-\boldsymbol{z}_{it}^{\top}\boldsymbol{\zeta^{i}}\!\leq s_{i},
\end{equation}
which is equivalent to $(8\text{b})$.
\qed
\noindent $\textbf{Remark 1}$
\;\emph{We note that since $f_{t}\left(\boldsymbol{x},\boldsymbol{\xi}\right)$ is convex in $\boldsymbol{\xi}$, then each function $G_{f_{t}}\left(\boldsymbol{z}_{it},\eta_{it},\boldsymbol{x}\right)$ in constraints $(8\text{b})$ is equivalent to maximizing a concave objective function on a closed convex set $\mathrm{\Xi}$.}\\
\noindent $\textbf{Remark 2}$
\;\emph{The exact reformulation $(8)$ of the set $\mathit{Z}_{D}$ is not convex because the index set $T(\boldsymbol{x})$ depends on $\boldsymbol{x}$ and each function $G_{f_{t}}\left(\boldsymbol{z}_{it},\eta_{it},\boldsymbol{x}\right)$ in constraints $(8\text{b})$ is not convex in general. Therefore, we attempt to investigate the tractability of the set $\mathit{Z}_{D}$ by establishing conditions under which $T(\boldsymbol{x})$ can be replaced by $\left[T\right]$ and $G_{f_{t}}\left(\boldsymbol{z}_{it},\eta_{it},\boldsymbol{x}\right)$ can be convex for any $i\in\left[N\right]$ and $t\in\left[T\right]$.}
\section{The Worst-Case CVaR Approximation}
\indent We first recall the definition of CVaR due to [22]. For a given measurable function $L(\boldsymbol{\xi}):\mathbb{R}^{m}\rightarrow\mathbb{R}$, let $\mathbb{P}$ be its probability distribution and the risk level $\epsilon\in\left(0,1\right)$, then the CVaR at level $\epsilon$ with respect to $\mathbb{P}$ is defined as
\begin{equation}\notag
\text{CVaR}_{1-\epsilon}\left(L\left(\boldsymbol{\xi}\right)\right):=\inf_{\beta\in\mathbb{R}}\left\{\beta+\frac{1}{\epsilon}\mathbb{E}_{\mathbb{P}}\left[\left(L\left(\boldsymbol{\xi}\right)-\beta\right)_{+}\right] \right\}.
\end{equation}
\indent It is well-known that for any $\boldsymbol{\alpha}\in \mathbb{R}_{++}^{T}$, the distributionally robust joint chance constraint $(1\text{c})$ can be reformulated as
\begin{equation}\tag{12}
\mathit{Z}_{D}=\left\{\boldsymbol{x}\in \mathbb{R}^n:\inf_{\mathbb{P}\in\mathcal{P}}\mathbb{P}\left\{\boldsymbol{\xi}:\max_{t\in\left[T\right]}\left\{\alpha_{t}(-f_{t}\left(\boldsymbol{x},\boldsymbol{\xi}\right))\right\}\leq 0\right\}\geq 1-\epsilon\right\}.
\end{equation}
\indent Note that $(12)$ represents a single distributionally robust chance constraint, which can be conservatively approximated by a worst-case CVaR constraint due to [7, 9]. Therefore, for any $\boldsymbol{\alpha}\in \mathbb{R}_{++}^{T}$, we have
\begin{equation}\notag
\mathit{Z}_{D}\supseteq\mathit{Z_{C}}(\boldsymbol{\alpha})=\left\{\boldsymbol{x}\in \mathbb{R}^n:\sup_{\mathbb{P}\in\mathcal{P}}\left\{\inf_{\beta\in\mathbb{R}}\left\{\beta+\frac{1}{\epsilon}\mathbb{E}_{\mathbb{P}}\left[\left(\max_{t\in\left[T\right]}\left\{\alpha_{t}\left(-f_{t}\left(\boldsymbol{x},\boldsymbol{\xi}\right)\right)\right\}-\beta\right)_{+}\right] \right\}\right\}\leq0\right\}.
\end{equation}
\indent Furthermore, we denote $\mathit{Z}_{C}$ as
\begin{equation}\notag
\mathit{Z}_{C}\!=\!\!\!\!\!\bigcup_{\boldsymbol{\alpha}\in \mathbb{R}_{++}^{T}}\!\!\!\!\!\mathit{Z_{C}}(\boldsymbol{\alpha})\!=\!\left\{\!\boldsymbol{x}\in \mathbb{R}^n\!:\!\boldsymbol{\alpha}\in \mathbb{R}_{++}^{T},\sup_{\mathbb{P}\in\mathcal{P}}\left\{\inf_{\beta\in\mathbb{R}}\left\{\!
\beta+\frac{1}{\epsilon}\mathbb{E}_{\mathbb{P}}\!\!\left[\left(\max_{t\in\left[T\right]}\left\{\alpha_{t}\left(-f_{t}\left(\boldsymbol{x},\boldsymbol{\xi}\right)\right)\right\}-\beta\right)_{+}\right] \right\}\right\}\leq0\!\right\}.
\end{equation}
\indent It should be noted that $\mathit{Z}_{C}=\bigcup_{\boldsymbol{\alpha}\in \mathbb{R}_{++}^{T}}\!\!\!\mathit{Z_{C}}(\boldsymbol{\alpha})\subseteq\mathit{Z_{D}}$, and then in the subsequent conclusion, we show that $\mathit{Z}_{C}$ can be reformulated as a disjunction of two sets $\mathit{Z}_{{C}_{1}}$ and $\mathit{Z}_{{C}_{2}}$.\\
\indent We first review the stochastic min-max theorem due to [23] before presenting the main results.
\begin{lemma}\label{66}
(Stochastic Min-max Equality, Theorem 2.1 in [23]) Let $\mathcal{A}$ be a nonempty (not necessarily convex) set of probability measures on measurable space $(\mathrm{\Xi}, \mathcal{B}(\mathrm{\Xi}))$ where $\mathrm{\Xi}\subseteq\mathbb{R}^m$ and $\mathcal{B}(\mathrm{\Xi})$ is the Borel $\sigma$-algebra. Assume that $\mathcal{A}$ is weakly compact. Let $\mathit{T}\subseteq\mathbb{R}^n$ be a closed convex set. Consider a function $\phi$: $\mathbb{R}^n\times \mathrm{\Xi}\rightarrow\mathbb{R}$. Assume that there exists a convex neighborhood $\mathit{V}$ of $\mathit{T}$ such that
for any $\boldsymbol{t}\in\mathit{V}$, the function $\phi(\boldsymbol{t},\cdot)$ is measurable, integrable with respect to all $\mathbb{P}\in\mathcal{A}$, and $\sup_{\mathbb{P}\in\mathcal{A}}\mathbb{E}_{\mathbb{P}}\left[\phi(\boldsymbol{t},\boldsymbol{\xi})\right]<\infty$. Further assume that $\phi(\cdot,\boldsymbol{\xi})$ is convex on $\mathit{V}$ for any $\boldsymbol{\xi}\in\mathrm{\Xi}$. Let $\bar{\boldsymbol{t}}\in\arg\min_{\boldsymbol{t}\in\mathit{T}}\sup_{\mathbb{P}\in\mathcal{A}}\mathbb{E}_{\mathbb{P}}\left[\phi(\boldsymbol{t},\boldsymbol{\xi})\right]$. Assume that for every $\boldsymbol{t}$ in a neighborhood of $\bar{\boldsymbol{t}}$, the function $\phi(\boldsymbol{t},\cdot)$ is bounded and upper-semicontinuous on $\mathrm{\Xi}$ and the function $\phi(\bar{\boldsymbol{t}},\cdot)$ is bounded and continuous on $\mathrm{\Xi}$. Then,
$$\inf_{\boldsymbol{t}\in\mathit{T}}\sup_{\mathbb{P}\in\mathcal{A}}\mathbb{E}_{\mathbb{P}}\left[\phi(\boldsymbol{t},\boldsymbol{\xi})\right]=\sup_{\mathbb{P}\in\mathcal{A}}\inf_{\boldsymbol{t}\in\mathit{T}}\mathbb{E}_{\mathbb{P}}\left[\phi(\boldsymbol{t},\boldsymbol{\xi})\right].$$
\end{lemma}

\indent We observe that Lemma 3 requires the ambiguity set to be weakly compact. This is indeed the case for Wasserstein ambiguity sets constructed from data due to [24].
\begin{lemma}\label{188}
(Proposition 3 in [24]) The Wasserstein ambiguity set $\mathcal{P}_{\mathit{W}}$ defined in $(5)$ is weakly compact.
\end{lemma}
\begin{proposition}\label{5}
The set $\mathit{Z}_{C}=\mathit{Z}_{{C}_{1}}\bigcup\mathit{Z}_{{C}_{2}}$, where
\begin{equation}\tag{13}
\mathit{Z}_{{C}_{1}}=\left\{\boldsymbol{x} \in \mathbb{R}^n:f_{t}\left(\boldsymbol{x},\boldsymbol{\xi}\right)\geq0, \forall \boldsymbol{\xi}\in\mathrm{\Xi}, \forall t\in\left[T\right]\right\},
\end{equation}
and
\begin{equation}\notag
~~~~~~~~~~~~~~~~~~~\mathit{Z}_{{C}_{2}}=
\left\{{\begin{array}{*{20}{c}}{\boldsymbol{x} \in \mathbb{R}^n:}&\begin{array}{l}
\lambda\delta+\frac{1}{N}\sum_{i=1}^{N}{s_{i}}\leq\epsilon,\\
G_{f_{t}}\left(\boldsymbol{z}_{it},\alpha_{t},\boldsymbol{x}\right)+1-\boldsymbol{z}_{it}^{\top}\boldsymbol{\zeta^{i}}-s_{i}\leq 0 ,~\forall i\in\left[N\right],\forall t\in\left[T\right],\\
\left \|\boldsymbol{z}_{it}\right \|_{*}-\lambda\leq0 ,~\forall i\in\left[N\right],\forall t\in\left[T\right],\\
\lambda\geq 0,~\alpha_{t}\geq 0,~s_{i}\geq0,~\forall i\in\left[N\right],\forall t\in\left[T\right],
\end{array}
\end{array}} \right\}\begin{array}{*{20}{c}}
{~~~~~~~~~~~\left( 14a \right)}\\
{~~~~~~~~~~~\left( 14b \right)}\\
{~~~~~~~~~~~\left( 14c \right)}\\
{~~~~~~~~~~~\left( 14d \right)}
\end{array}
\end{equation}
where $G_{f_{t}}\left(\boldsymbol{z}_{it},\alpha_{t},\boldsymbol{x}\right)=
\sup_{\boldsymbol{\xi}\in\mathrm{\Xi}}\left[\boldsymbol{z}_{it}^{\top}\boldsymbol{\xi}
-\alpha_{t}f_{t}\left(\boldsymbol{x},\boldsymbol{\xi}\right)\right]$.
\end{proposition}
\noindent{\it Proof}.~\;We separate the proof into three parts.\\
\noindent$(i)$ Note that
\begin{alignat}{2}
\mathit{Z}_{C}&=\left\{\boldsymbol{x}\in \mathbb{R}^n:\boldsymbol{\alpha}\in \mathbb{R}_{++}^{T},\sup_{\mathbb{P}\in\mathcal{P}}\left\{\inf_{\beta\in\mathbb{R}}\left\{\beta+\frac{1}{\epsilon}\mathbb{E}_{\mathbb{P}}\left[\left(\max_{t\in\left[T\right]}\left\{\alpha_{t}\left(-f_{t}\left(\boldsymbol{x},\boldsymbol{\xi}\right)\right)\right\}-\beta\right)_{+}\right] \right\}\right\}\leq0\right\}&\tag{15}\\
&=\left\{\boldsymbol{x}\in \mathbb{R}^n:\boldsymbol{\alpha}\in \mathbb{R}_{++}^{T},\inf_{\beta\in\mathbb{R}}\left\{\beta+\frac{1}{\epsilon}\sup_{\mathbb{P}\in\mathcal{P}}\mathbb{E}_{\mathbb{P}}\left[\left(\max_{t\in\left[T\right]}\left\{\alpha_{t}\left(-f_{t}\left(\boldsymbol{x},\boldsymbol{\xi}\right)\right)\right\}-\beta\right)_{+}\right] \right\}\leq0\right\}&\tag{16}\\
&=\left\{\boldsymbol{x}\in \mathbb{R}^n:\boldsymbol{\alpha}\in \mathbb{R}_{++}^{T},\beta\in \mathbb{R},\beta+\frac{1}{\epsilon}\sup_{\mathbb{P}\in\mathcal{P}}\mathbb{E}_{\mathbb{P}}\left[\left(\max_{t\in\left[T\right]}\left\{\alpha_{t}\left(-f_{t}\left(\boldsymbol{x},\boldsymbol{\xi}\right)\right)\right\}-\beta\right)_{+}\right] \leq0\right\},&\tag{17}
\end{alignat}
where the second equality is due to Lemma 3 and Lemma 4 and the third equality follows from replacing infimum operator with its equivalent ``existence'' argument.\\
\indent Then, we prove that $\beta\leq0$. Suppose that $\beta>0$. Since $\left(\max_{t\in\left[T\right]}\left\{\alpha_{t}\left(-f_{t}\left(\boldsymbol{x},\boldsymbol{\xi}\right)\right)\right\}-\beta\right)_{+}\geq0$ for any $\boldsymbol{\xi}\in\mathrm{\Xi}$, we must have $\mathbb{E}_{\mathbb{P}}\left[\left(\max_{t\in\left[T\right]}\left\{\alpha_{t}\left(-f_{t}\left(\boldsymbol{x},\boldsymbol{\xi}\right)\right)\right\}-\beta\right)_{+}\right] \geq0$. Thus, the left-hand side of $(17)$ is strictly positive, which yields a contradiction.\\
\noindent$(ii)$ Now, we show that $\mathit{Z}_{C}\subseteq \mathit{Z}_{{C}_{1}}\bigcup \mathit{Z}_{{C}_{2}}$. For any $\boldsymbol{x}\in \mathit{Z}_{C}$, there exists
$(\boldsymbol{\alpha},\beta)\in\mathbb{R}_{++}^{T}\times\mathbb{R}_{-}$ such that
\begin{equation}\tag{18}
\beta\epsilon+\sup_{\mathbb{P}\in\mathcal{P}}\mathbb{E}_{\mathbb{P}}\left[\left(\max_{t\in\left[T\right]}\left\{\alpha_{t}\left(-f_{t}\left(\boldsymbol{x},\boldsymbol{\xi}\right)\right)\right\}-\beta\right)_{+}\right] \leq0.
\end{equation}
\indent Then, we distinguish whether $\beta=0$ or $\beta<0$.\\
\indent Case 1.
Note that if $\beta=0$, then inequality $(18)$ implies
\begin{equation}\tag{19}
\sup_{\mathbb{P}\in\mathcal{P}}\mathbb{E}_{\mathbb{P}}\left[\left(\max_{t\in\left[T\right]}\left\{\alpha_{t}\left(-f_{t}\left(\boldsymbol{x},\boldsymbol{\xi}\right)\right)\right\}\right)_{+}\right] \leq0,
\end{equation}
which is equivalent to
\begin{equation}\notag
\inf_{\mathbb{P}\in\mathcal{P}}\mathbb{P}\left[f_{t}\left(\boldsymbol{x},\boldsymbol{\xi}\right)\geq0,\forall t\in\left[T\right]\right]=1>1-\epsilon,
\end{equation}
and hence by continuity of each function $f_{t}\left(\boldsymbol{x},\boldsymbol{\xi}\right)$, we have $f_{t}\left(\boldsymbol{x},\boldsymbol{\xi}\right)\geq0$ for any $\boldsymbol{\xi}\in\mathrm{\Xi}$. Thus,
$\boldsymbol{x}\in\mathit{Z}_{{C}_{1}}$.\\
\indent Case 2.
On the other hand, if $\beta<0$, then divide $(18)$ by $-\beta$ and add $\epsilon$ on both sides, we have
\begin{equation}\tag{20}
-\frac{1}{\beta}\sup_{\mathbb{P}\in\mathcal{P}}\mathbb{E}_{\mathbb{P}}\left[\left(\max_{t\in\left[T\right]}\left\{\alpha_{t}\left(-f_{t}\left(\boldsymbol{x},\boldsymbol{\xi}\right)\right)\right\}-\beta\right)_{+}\right] \leq\epsilon.
\end{equation}
\indent Since $\beta<0$, we can redefine $\alpha_{t}$ as $\alpha_{t}/(-\beta)$ for any $t\in\left[T\right]$. Then, inequality $(20)$ can be rewritten as
\begin{equation}\tag{21}
\sup_{\mathbb{P}\in\mathcal{P}}\mathbb{E}_{\mathbb{P}}\left[\left(\max_{t\in\left[T\right]}\left\{\alpha_{t}\left(-f_{t}\left(\boldsymbol{x},\boldsymbol{\xi}\right)\right)\right\}+1\right)_{+}\right] \leq\epsilon.
\end{equation}
\indent Since
\begin{equation}\notag
\begin{aligned}
~~~~~~~~~~~\left(\max_{t\in\left[T\right]}\left\{\alpha_{t}\left(-f_{t}\left(\boldsymbol{x},\boldsymbol{\xi}\right)\right)\right\}+1\right)_{+}
&=\max\left\{\max_{t\in\left[T\right]}\left\{\alpha_{t}
\left(-f_{t}\left(\boldsymbol{x},\boldsymbol{\xi}\right)\right)\right\}+1,\,0\right\}\\
&=\max\left\{\alpha_{1}
\left(-f_{1}\left(\boldsymbol{x},\boldsymbol{\xi}\right)\right)+1,\cdots,\alpha_{T}
\left(-f_{T}\left(\boldsymbol{x},\boldsymbol{\xi}\right)\right)+1,\,0\right\}.
\end{aligned}
\end{equation}
\indent Therefore, by Lemma 1, for any given $\boldsymbol{\alpha}\in \mathbb{R}_{++}^{T}$, the supremum in the left-hand side of $(21)$ is equivalent to
\begin{alignat}{2}
~\inf_{\lambda,\boldsymbol{s},\boldsymbol{z},\boldsymbol{\alpha}}~&\lambda\delta+\frac{1}{N}\sum_{i=1}^{N}{s_{i}}& ~\tag{22a}\\
\mbox{s.t.}\quad
&\sup_{\boldsymbol{\xi}\in\mathrm{\Xi}}\left[\boldsymbol{z}_{it}^{\top}\boldsymbol{\xi}-\alpha_{t}f_{t}\left(\boldsymbol{x},\boldsymbol{\xi}\right)+1\right]
-\boldsymbol{z}_{it}^{\top}\boldsymbol{\zeta^{i}}\leq s_{i},~\forall i\in\left[N\right],\forall t\in\left[T\right], & ~\tag{22b}\\
&\left \|\boldsymbol{z}_{it}\right \|_{*}\leq \lambda,~\forall i\in\left[N\right],\forall t\in\left[T\right],&~\tag{22c}\\
&\lambda\!\geq 0,~\alpha_{t}\!>0,~s_{i}\geq0,~\forall i\in\left[N\right],\forall t\in\left[T\right],&~\tag{22d}
\end{alignat}
which implies that $\boldsymbol{x}\in\mathit{Z}_{{C}_{2}}$.\\
\noindent$(iii)$ Now, we show that $\mathit{Z}_{{C}_{1}}\bigcup \mathit{Z}_{{C}_{2}}\subseteq \mathit{Z}_{C}$. Similarly, given $\boldsymbol{x}\in \mathit{Z}_{{C}_{1}}\bigcup \mathit{Z}_{{C}_{2}}$. If $\boldsymbol{x}\in \mathit{Z}_{{C}_{1}}$, then we let $\beta=0$, $\boldsymbol{\alpha}=\boldsymbol{e}$, thus $\boldsymbol{x}\in \mathit{Z}_{C}$.  If $\boldsymbol{x}\in\mathit{Z}_{{C}_{2}}$, there exists $\left(\lambda',\boldsymbol{s}',\boldsymbol{z}',\boldsymbol{\alpha}',\boldsymbol{x}\right)$ which satisfies the constraints in $(14)$. Since $\boldsymbol{\alpha}'>\boldsymbol{0}$, which could be confirmed in the following Corollary 1, and hence let $\beta=-1$, $\boldsymbol{\alpha}=\boldsymbol{\alpha}'$ in $(17)$. Then, by strong duality result introduced in Lemma 1, we have $\boldsymbol{x}\in \mathit{Z}_{C}$.
\qed
\noindent $\textbf{Remark 3}$
\;\emph{To solve the inner approximation of optimization problem $(1)$ $\left(\text{i.e}.,\min_{\boldsymbol{x}\in\mathit{S}\bigcap\mathit{Z}_{C}}\boldsymbol{c}^{\top}\boldsymbol{x}\right)$, we can optimize $\boldsymbol{c}^{\top}\boldsymbol{x}$ over $\mathit{S}\bigcap\mathit{Z}_{{C}_{1}}$ and $\mathit{S}\bigcap\mathit{Z}_{{C}_{2}}$ separately, then choose the smallest value.}\\
\noindent $\textbf{Remark 4}$
\;\emph{We observe that each function $G_{f_{t}}\left(\boldsymbol{z}_{it},\alpha_{t},\boldsymbol{x}\right)$ is merely biconvex, but not jointly convex in $\boldsymbol{x}$ and $\alpha_{t}$. Then, the left-hand sides of the constraint system $(14)$ are biconvex in $\boldsymbol{\alpha}$ and $\left(\lambda,\boldsymbol{s},\boldsymbol{z},\boldsymbol{x}\right)$, i.e., they are convex in $\left(\lambda,\boldsymbol{s},\boldsymbol{z},\boldsymbol{x}\right)$ for any given $\boldsymbol{\alpha}\in\mathbb{R}_{+}^{T}$, and also convex in $\boldsymbol{\alpha}$ for any given $\left(\lambda,\boldsymbol{s},\boldsymbol{z},\boldsymbol{x}\right)$. Thus, optimization problem $\min_{\boldsymbol{x}\in\mathit{S}\bigcap\mathit{Z}_{{C}_{2}}}\boldsymbol{c}^{\top}\boldsymbol{x}$ is non-convex.}\\
\indent Next, we prove that in the constraint system $(14)$, $\alpha_{t}$ must be strictly positive for any $t\in\left[T\right]$.
\begin{corollary}\label{3}
For any $\boldsymbol{x}$ satisfying $(14)$, we must have $\alpha_{t}>0$ for any $t\in\left[T\right]$.
\end{corollary}
{\it Proof}.~\; Suppose that we let $\alpha_{{t}_{0}}=0$ for some $t_{0}\in\left[T\right]$, then from $(14\text{b})$, we have
\begin{equation}\notag
G_{f_{{t}_{0}}}\left(\boldsymbol{z}_{it_{0}},0,\boldsymbol{x}\right)+1-\boldsymbol{z}_{it_{0}}^{\top}\boldsymbol{\zeta^{i}}=
\sup_{\boldsymbol{\xi}\in\mathrm{\Xi}}\boldsymbol{z}_{it_{0}}^{\top}\boldsymbol{\xi}+1-\boldsymbol{z}_{it_{0}}^{\top}\boldsymbol{\zeta^{i}}
=\sup_{\boldsymbol{\xi}\in\mathrm{\Xi}}\left[\boldsymbol{z}_{it_{0}}^{\top}\left(\boldsymbol{\xi}-\boldsymbol{\zeta^{i}}\right)\right]+1\leq s_{i},~\forall i\in\left[N\right].
\end{equation}
\indent We note that in $(14\text{c})$, $\left \|\boldsymbol{z}_{it_{0}}\right \|_{*}\leq \lambda$ for any $i\in\left[N\right]$, we may thus conclude that
\begin{equation}\notag
\inf_{\left \|\boldsymbol{z}_{it_{0}}\right \|_{*}\leq \lambda}\left[\boldsymbol{z}_{it_{0}}^{\top}\left(\boldsymbol{\xi}-\boldsymbol{\zeta^{i}}\right)\right]+1\leq s_{i},~\forall \boldsymbol{\xi}\in\mathrm{\Xi},\forall i\in\left[N\right],
\end{equation}
which can be rewritten as
\begin{equation}\notag
-\!\!\!\!\!\!\sup_{\left \|\boldsymbol{z}_{it_{0}}\right \|_{*}\leq \lambda}\left[-\boldsymbol{z}_{it_{0}}^{\top}\left(\boldsymbol{\xi}-\boldsymbol{\zeta^{i}}\right)\right]+1\leq s_{i},~\forall \boldsymbol{\xi}\in\mathrm{\Xi},\forall i\in\left[N\right],
\end{equation}
furthermore, which is equivalent to
\begin{equation}\notag
-\lambda\left\|\boldsymbol{\zeta^{i}}-\boldsymbol{\xi}\right\|+1\leq s_{i},~\forall \boldsymbol{\xi}\in\mathrm{\Xi}, \forall i\in\left[N\right].
\end{equation}
\indent Therefore, according to $(14\text{a})$, we have
\begin{equation}\notag
\lambda\delta-\frac{1}{N}\sum_{i=1}^{N}\lambda\left\|\boldsymbol{\zeta^{i}}-\boldsymbol{\xi}\right\|+1\leq\lambda\delta+\frac{1}{N}\sum_{i=1}^{N}{s_{i}}\leq\epsilon,~\forall \boldsymbol{\xi}\in\mathrm{\Xi},
\end{equation}
which can be rewritten as
\begin{equation}\tag{23}
\lambda\delta-\frac{1}{N}\sum_{i=1}^{N}\lambda\inf_{\boldsymbol{\xi}\in\mathrm{\Xi}}\left\|\boldsymbol{\zeta^{i}}-\boldsymbol{\xi}\right\|+1\leq\lambda\delta+\frac{1}{N}\sum_{i=1}^{N}{s_{i}}\leq\epsilon.
\end{equation}
\indent Note that $\inf_{\boldsymbol{\xi}\in\mathrm{\Xi}}\left\|\boldsymbol{\zeta^{i}}-\boldsymbol{\xi}\right\|=0$ for any $i\in\left[N\right]$, thus inequalities $(23)$ yields a contradiction to the fact that $\epsilon<1$.
\qed

\indent Meanwhile, we also show that $\alpha_{t}$ could be bounded for any $t\in\left[T\right]$.
\begin{corollary}\label{6}
If $\mathit{S}$ is compact and $\left[T\right]=T(\boldsymbol{x})$ for any $\boldsymbol{x}\in\mathit{Z}_{{C}_{2}}$, then there exists an $\boldsymbol{M}\in\mathbb{R}_{++}^{T}$ such that $\alpha_{t}\leq M_{t}$ for any $t\in\left[T\right]$.
\end{corollary}
\noindent{\it Proof}.~\; We note that the statement that $\left[T\right]=T(\boldsymbol{x})$ for any $\boldsymbol{x}\in\mathit{Z}_{{C}_{2}}$ implies given $\boldsymbol{x}\in\mathit{S}\bigcap\mathit{Z}_{{C}_{2}}$, for any $t\in\left[T\right]$, there exists $\boldsymbol{\xi}\in\mathrm{\Xi}$ such that $f_{t}(\boldsymbol{x},\boldsymbol{\xi})<0$.\\
\indent Since constraints $(14\text{b})$ can be rewritten as
\begin{equation}\notag
\sup_{\boldsymbol{\xi}\in\mathrm{\Xi}}\left[\boldsymbol{z}_{it}^{\top}\boldsymbol{\xi}
-\alpha_{t}f_{t}\left(\boldsymbol{x},\boldsymbol{\xi}\right)\right]\leq-\left(1-\boldsymbol{z}_{it}^{\top}\boldsymbol{\zeta}^{i}-s_{i}\right) ,~\forall i\in\left[N\right],\forall t\in\left[T\right],
\end{equation}
furthermore, which is equivalent to
\begin{equation}\notag
\boldsymbol{z}_{it}^{\top}\boldsymbol{\xi}
-\alpha_{t}f_{t}\left(\boldsymbol{x},\boldsymbol{\xi}\right)\leq-\left(1-\boldsymbol{z}_{it}^{\top}\boldsymbol{\zeta}^{i}-s_{i}\right) ,~\forall\boldsymbol{\xi}\in\mathrm{\Xi},~\forall i\in\left[N\right], \forall t\in\left[T\right].
\end{equation}
\indent Since $\mathit{S}$ is compact and $\mathrm{\Xi}$ is closed, we have $f_{t}(\boldsymbol{x},\boldsymbol{\xi})$ must be finite. Therefore, when $f_{t}(\boldsymbol{x},\boldsymbol{\xi})<0$, we obtain
\begin{equation}\notag
\alpha_{t}\leq\frac{1}{f_{t}\left(\boldsymbol{x},\boldsymbol{\xi}\right)}\left[\boldsymbol{z}_{it}^{\top}\boldsymbol{\xi} +1-\boldsymbol{z}_{it}^{\top}\boldsymbol{\zeta}^{i}-s_{i}\right] ,~\forall i\in\left[N\right],\forall t\in\left[T\right].
\end{equation}
\indent Thus, one can find an upper bound $\boldsymbol{M}\in\mathbb{R}_{++}^{T}$ such that $\alpha_{t}\leq M_{t}$ for any $t\in\left[T\right]$.
\qed
\indent The following theorem is derived by constructing two new decision variables which allows us to eliminate biconvex terms and proving that the new formulation is a convex relaxation of the set $\mathit{Z}_{{C}_{2}}$.
\begin{theorem}\label{10077}
The set $\mathit{Z}_{{C}_{2}}$ can be outer approximated by
\begin{equation}\notag
~~~~~~~~~~~~~~~~~~~~~\mathit{\tilde{Z}}_{{C}_{2}}=
\left\{{\begin{array}{*{20}{c}}{\boldsymbol{x} \in \mathbb{R}^n:}&\begin{array}{l}
\left \|\boldsymbol{v}_{it}\right \|_{*}\delta+\frac{1}{N}\sum_{i=1}^{N}{q_{it}}\leq\epsilon\alpha_{t},~\forall i\in\left[N\right],\forall t\in\left[T\right],\\
G_{f_{t}}\left(\boldsymbol{v}_{it},1,\boldsymbol{x}\right)+\alpha_{t}-\boldsymbol{v}_{it}^{\top}\boldsymbol{\zeta^{i}}-q_{it}\leq0,~\forall i\in\left[N\right],\forall t\in\left[T\right],\\
\alpha_{t}> 0,~q_{it}\geq0,~\forall i\in\left[N\right],\forall t\in\left[T\right],
\end{array}
\end{array}}\!\!\right\}\begin{array}{*{20}{c}}
{~~~~~~~~~~\left( 24a \right)}\\
{~~~~~~~~~~\left( 24b \right)}\\
{~~~~~~~~~~\left( 24c \right)}
\end{array}
\end{equation}
which is a convex set.
\end{theorem}
\noindent{\it Proof}.~\; By Proposition 1 and Corollary 1, the set $\mathit{Z}_{{C}_{2}}$ can be expressed as
\begin{equation}\notag
\mathit{Z}_{{C}_{2}}\!=\!
\left\{\!\!{\begin{array}{*{20}{c}}{\boldsymbol{x} \in \mathbb{R}^n\!:\!}&\begin{array}{l}
\lambda\delta+\frac{1}{N}\sum_{i=1}^{N}{s_{i}}\leq\epsilon,\\
\alpha_{t}\sup_{\boldsymbol{\xi}\in\mathrm{\Xi}}\left[(\frac{\boldsymbol{z}_{it}}{\alpha_{t}})^{\top}\boldsymbol{\xi}
-f_{t}\left(\boldsymbol{x},\boldsymbol{\xi}\right)\right]+1-\alpha_{t}(\frac{\boldsymbol{z}_{it}}{\alpha_{t}})^{\top}\boldsymbol{\zeta^{i}}-s_{i}\leq 0 ,~\forall i\in\left[N\right],\forall t\in\left[T\right],\\
\left \|\boldsymbol{z}_{it}\right \|_{*}-\lambda\leq0 ,~\forall i\in\left[N\right],\forall t\in\left[T\right],\\
\lambda\geq 0,~\alpha_{t}>0,~s_{i}\geq0,~\forall i\in\left[N\right],\forall t\in\left[T\right].
\end{array}
\end{array}\!\!} \right\}
\end{equation}
\indent Note that for any $\boldsymbol{x}\in\mathit{Z}_{{C}_{2}}$, there exists $\left(\lambda,\boldsymbol{s},\boldsymbol{z},\boldsymbol{\alpha},\boldsymbol{x}\right)$ which satisfies the constraints above, and then, we could let $\lambda=\max_{i\in\left[N\right],t\in\left[T\right]}\left\|\boldsymbol{z}_{it}\right \|_{*}$, which finally allows us to reformulate $\mathit{Z}_{{C}_{2}}$ as
\begin{equation}\notag
\mathit{Z}_{{C}_{2}}=
\left\{{\begin{array}{*{20}{c}}{\boldsymbol{x} \in \mathbb{R}^n:}&\begin{array}{l}
\left \|\boldsymbol{z}_{it}\right \|_{*}\delta+\frac{1}{N}\sum_{i=1}^{N}{s_{i}}\leq\epsilon,~\forall i\in\left[N\right],\forall t\in\left[T\right],\\
\sup_{\boldsymbol{\xi}\in\mathrm{\Xi}}\left[(\frac{\boldsymbol{z}_{it}}{\alpha_{t}})^{\top}\boldsymbol{\xi}
-f_{t}\left(\boldsymbol{x},\boldsymbol{\xi}\right)\right]+\frac{1}{\alpha_{t}}-(\frac{\boldsymbol{z}_{it}}{\alpha_{t}})^{\top}\boldsymbol{\zeta^{i}}-\frac{s_{i}}{\alpha_{t}}\leq 0 ,~\forall i\in\left[N\right],\forall t\in\left[T\right],\\
\alpha_{t}>0,~s_{i}\geq0,~\forall i\in\left[N\right],\forall t\in\left[T\right].
\end{array}
\end{array}} \right\}
\end{equation}
\indent Subsequently, we introduce two new decision variables $\boldsymbol{v}_{it}=\frac{\boldsymbol{z}_{it}}{\alpha_{t}}$ and $q_{it}=\frac{s_{i}}{\alpha_{t}}$. Then, the set $\mathit{Z}_{{C}_{2}}$ can be outer approximated by
\begin{equation}\notag
\mathit{\tilde{Z}}_{{C}_{2}}=
\left\{{\begin{array}{*{20}{c}}{\boldsymbol{x} \in \mathbb{R}^n:}&\begin{array}{l}
\left \|\boldsymbol{v}_{it}\right \|_{*}\delta+\frac{1}{N}\sum_{i=1}^{N}{q_{it}}\leq\frac{\epsilon}{\alpha_{t}},~\forall i\in\left[N\right],\forall t\in\left[T\right],\\
\sup_{\boldsymbol{\xi}\in\mathrm{\Xi}}\left[\boldsymbol{v}_{it}^{\top}\boldsymbol{\xi}
-f_{t}\left(\boldsymbol{x},\boldsymbol{\xi}\right)\right]+\frac{1}{\alpha_{t}}-\boldsymbol{v}_{it}^{\top}\boldsymbol{\zeta^{i}}-q_{it}\leq 0 ,~\forall i\in\left[N\right],\forall t\in\left[T\right],\\
\alpha_{t}>0,~q_{it}\geq0,~\forall i\in\left[N\right],\forall t\in\left[T\right].
\end{array}
\end{array}} \right\}
\end{equation}
\indent We now perform variable substitution in which we replace $\frac{1}{\alpha_{t}}$ by $\alpha_{t}$, which yields the following reformulation of $\mathit{\tilde{Z}}_{{C}_{2}}$
\begin{equation}\notag
\mathit{\tilde{Z}}_{{C}_{2}}=
\left\{{\begin{array}{*{20}{c}}{\boldsymbol{x} \in \mathbb{R}^n:}&\begin{array}{l}
\left \|\boldsymbol{v}_{it}\right \|_{*}\delta+\frac{1}{N}\sum_{i=1}^{N}{q_{it}}\leq\epsilon\alpha_{t},~\forall i\in\left[N\right],\forall t\in\left[T\right],\\
\sup_{\boldsymbol{\xi}\in\mathrm{\Xi}}\left[\boldsymbol{v}_{it}^{\top}\boldsymbol{\xi}
-f_{t}\left(\boldsymbol{x},\boldsymbol{\xi}\right)\right]+\alpha_{t}-\boldsymbol{v}_{it}^{\top}\boldsymbol{\zeta^{i}}-q_{it}\leq 0 ,~\forall i\in\left[N\right],\forall t\in\left[T\right],\\
\alpha_{t}>0,~q_{it}\geq0,~\forall i\in\left[N\right],\forall t\in\left[T\right],
\end{array}
\end{array}} \right\}
\end{equation}
and thus the claim follows.
\qed
\noindent $\textbf{Remark 5}$
\;\emph{It is easy to check that if there is a single uncertain constraint $\left(\text{i.e}., T=1\right)$, then $\mathit{\tilde{Z}}_{{C}_{2}}=\mathit{Z}_{{C}_{2}}$.}\\
\indent We now present tractable conic reformulations of the set $\mathit{\tilde{Z}}_{{C}_{2}}$
when the constraint function is quadratic convex in the uncertainty and concave in the decision variable,
and the support of the uncertainty is polyhedron or ellipsoid. In addition, we provide the detailed proofs in Appendix A.
\begin{proposition}
If the support set $\mathrm{\Xi}$ is polyhedral, i.e., $\mathrm{\Xi}=\{\boldsymbol{\xi}\in\mathbb{R}^{m}:\boldsymbol{a}_{k}^{\top}\boldsymbol{\xi}\leq d_{k}, \boldsymbol{a}_{k}\in\mathbb{R}^{m}, d_{k}>0, k=1,2,\cdots,l\}$, and $f_{t}\left(\boldsymbol{x},\boldsymbol{\xi}\right)=\boldsymbol{\xi}^{\top}\boldsymbol{x}+
\langle\boldsymbol{A}^{t},\boldsymbol{\xi}\boldsymbol{\xi}^{\top}\rangle+(\boldsymbol{b}^{t})^{\top}\boldsymbol{x}+h^{t}$ for any $t\in \left[T\right]$, where $\boldsymbol{A}^{t}\succeq0$. Then
\begin{equation}\notag
~~~~~~~~~~~~~~~\mathit{\tilde{Z}}_{{C}_{2}}=
\left\{{\begin{array}{*{20}{c}}{\boldsymbol{x} \in \mathbb{R}^n:}&\begin{array}{l}
\left \|\boldsymbol{v}_{it}\right \|_{*}\delta+\frac{1}{N}\sum_{i=1}^{N}{q_{it}}\leq\epsilon\alpha_{t},~\forall i\in\left[N\right],\forall t\in\left[T\right],\\
u_{it}-\left[(\boldsymbol{b}^{t})^{\top}\boldsymbol{x}+h^{t}\right]+\alpha_{t}-\boldsymbol{v}_{it}^{\top}\boldsymbol{\zeta^{i}}\leq q_{it},~\forall i\in\left[N\right],\forall t\in\left[T\right],\\
U(u_{it},\boldsymbol{v}_{it},\boldsymbol{x},\boldsymbol{\nu}^{it})\succeq0,~\forall i\in\left[N\right],\forall t\in\left[T\right],\\
\alpha_{t}> 0,~\nu^{it}_{k}\geq0,~q_{it}\geq0,u_{it}\in\mathbb{R},~\forall i\in\left[N\right],\forall t\in\left[T\right],\forall k\in\left[l\right],
\end{array}
\end{array}}\!\!\right\}\begin{array}{*{20}{c}}
{~~~~~~~~\left( 25a \right)}\\
{~~~~~~~~\left( 25b \right)}\\
{~~~~~~~~\left( 25c \right)}\\
{~~~~~~~~\left( 25d \right)}
\end{array}
\end{equation}
where $U(u_{it},\boldsymbol{v}_{it},\boldsymbol{x},\boldsymbol{\nu}^{it})=
\begin{bmatrix}
\boldsymbol{A}^{t}&-\frac{1}{2}(\boldsymbol{v}_{it}-\boldsymbol{x}-\sum_{k=1}^{l}\nu_{k}^{it}\boldsymbol{a}_{k})\\
-\frac{1}{2}(\boldsymbol{v}_{it}-\boldsymbol{x}-\sum_{k=1}^{l}\nu_{k}^{it}\boldsymbol{a}_{k})^{\top}&u_{it}-\sum_{k=1}^{l}\nu_{k}^{it}d_{k}
\end{bmatrix}$.
\end{proposition}
\begin{proposition}
If the support set $\mathrm{\Xi}$ is ellipsoidal, i.e., $\mathrm{\Xi}=\{\boldsymbol{\xi}\in\mathbb{R}^{m}:({\boldsymbol{\xi}-\boldsymbol{\xi}_{0}})^{\top}\boldsymbol{W}^{-1}({\boldsymbol{\xi}-\boldsymbol{\xi}_{0}})\leq 1\}$, where $\boldsymbol{W}\succ0$, and $f_{t}\left(\boldsymbol{x},\boldsymbol{\xi}\right)=\boldsymbol{\xi}^{\top}\boldsymbol{x}+\langle\boldsymbol{A}^{t},\boldsymbol{\xi}\boldsymbol{\xi}^{\top}\rangle+(\boldsymbol{b}^{t})^{\top}\boldsymbol{x}+h^{t}$ for any $t\in \left[T\right]$, where $\boldsymbol{A}^{t}\succeq0$. Then
\begin{equation}\notag
~~~~~~~~~~~~~~~~~\mathit{\tilde{Z}}_{{C}_{2}}=
\left\{{\begin{array}{*{20}{c}}{\boldsymbol{x} \in \mathbb{R}^n:}&\begin{array}{l}
\left \|\boldsymbol{v}_{it}\right \|_{*}\delta+\frac{1}{N}\sum_{i=1}^{N}{q_{it}}\leq\epsilon\alpha_{t},~\forall i\in\left[N\right],\forall t\in\left[T\right],\\
u_{it}-\left[(\boldsymbol{b}^{t})^{\top}\boldsymbol{x}+h^{t}\right]+\alpha_{t}-\boldsymbol{v}_{it}^{\top}\boldsymbol{\zeta^{i}}\leq q_{it},~\forall i\in\left[N\right],\forall t\in\left[T\right],\\
U(u_{it},\boldsymbol{v}_{it},\boldsymbol{x},\nu_{it})\succeq0,~\forall i\in\left[N\right],\forall t\in\left[T\right],\\
\alpha_{t}> 0,~\nu_{it}\geq0,~q_{it}\geq0,u_{it}\in\mathbb{R},~\forall i\in\left[N\right],\forall t\in\left[T\right],
\end{array}
\end{array}}\!\!\right\}\begin{array}{*{20}{c}}
{~~~~~~\left( 26a \right)}\\
{~~~~~~\left( 26b \right)}\\
{~~~~~~\left( 26c \right)}\\
{~~~~~~\left( 26d \right)}
\end{array}
\end{equation}
where $U(u_{it},\boldsymbol{v}_{it},\boldsymbol{x},\nu_{it})=
\begin{bmatrix}
\boldsymbol{A}^{t}+\nu_{it}\boldsymbol{W}^{-1}&-\frac{1}{2}(2\nu_{it}\boldsymbol{W}^{-1}\boldsymbol{\xi}_{0}+\boldsymbol{v}_{it}-\boldsymbol{x})\\
-\frac{1}{2}(2\nu_{it}\boldsymbol{W}^{-1}\boldsymbol{\xi}_{0}+\boldsymbol{v}_{it}-\boldsymbol{x})^{\top}&u_{it}+\nu_{it}\boldsymbol{\xi}_{0}^{\top}\boldsymbol{W}^{-1}\boldsymbol{\xi}_{0}-\nu_{it}
\end{bmatrix}$.
\end{proposition}
\begin{proposition}
If the support set $\mathrm{\Xi}$ is ellipsoidal, i.e., $\mathrm{\Xi}=\{\boldsymbol{\xi}\in\mathbb{R}^{m}:({\boldsymbol{\xi}-\boldsymbol{\xi}_{0}})^{\top}\boldsymbol{W}^{-1}({\boldsymbol{\xi}-\boldsymbol{\xi}_{0}})\leq 1\}$, where $\boldsymbol{W}\succ0$, and $f_{t}\left(\boldsymbol{x},\boldsymbol{\xi}\right)=\boldsymbol{w}_{t}(\boldsymbol{\xi})^{\top}\boldsymbol{x}$ for any $t\in \left[T\right]$, where each component $w_{tj}(\boldsymbol{\xi})$ of $\boldsymbol{w}_{t}(\boldsymbol{\xi})$ is quadratic convex in $\boldsymbol{\xi}$, i.e., it has the form $w_{tj}(\boldsymbol{\xi})=\boldsymbol{\xi}^{\top}\boldsymbol{W}_{tj}\boldsymbol{\xi}+\boldsymbol{r}_{tj}^{\top}\boldsymbol{\xi}+h_{tj}$, where $\boldsymbol{W}_{tj}\succeq0$. Then
\begin{equation}\notag
~~~~~~~~~~~~~~~~~~~~~\mathit{\tilde{Z}}_{{C}_{2}}=
\left\{{\begin{array}{*{20}{c}}{\boldsymbol{x} \in \mathbb{R}^n:}&\begin{array}{l}
\left \|\boldsymbol{v}_{it}\right \|_{*}\delta+\frac{1}{N}\sum_{i=1}^{N}{q_{it}}\leq\epsilon\alpha_{t},~\forall i\in\left[N\right],\forall t\in\left[T\right],\\
u_{it}-H_{t}(\boldsymbol{x})+\alpha_{t}-\boldsymbol{v}_{it}^{\top}\boldsymbol{\zeta^{i}}\leq q_{it},~\forall i\in\left[N\right],\forall t\in\left[T\right],\\
U(u_{it},\boldsymbol{v}_{it},\boldsymbol{x},\nu_{it})\succeq0,~\forall i\in\left[N\right],\forall t\in\left[T\right],\\
\alpha_{t}> 0,~\nu_{it}\geq0,~q_{it}\geq0,u_{it}\in\mathbb{R},~\forall i\in\left[N\right],\forall t\in\left[T\right],
\end{array}
\end{array}}\!\!\right\}\begin{array}{*{20}{c}}
{~~~~~~~~\left( 27a \right)}\\
{~~~~~~~~\left( 27b \right)}\\
{~~~~~~~~\left( 27c \right)}\\
{~~~~~~~~\left( 27d \right)}
\end{array}
\end{equation}
where $$U(u_{it},\boldsymbol{v}_{it},\boldsymbol{x},\nu_{it})=
\begin{bmatrix}
\boldsymbol{W}_{t}(\boldsymbol{x})+\nu_{it}\boldsymbol{W}^{-1}&-\frac{1}{2}(2\nu_{it}\boldsymbol{W}^{-1}\boldsymbol{\xi}_{0}+\boldsymbol{v}_{it}-\boldsymbol{R}_{t}(\boldsymbol{x}))\\
-\frac{1}{2}(2\nu_{it}\boldsymbol{W}^{-1}\boldsymbol{\xi}_{0}+\boldsymbol{v}_{it}-\boldsymbol{R}_{t}(\boldsymbol{x}))^{\top}&u_{it}+\nu_{it}\boldsymbol{\xi}_{0}^{\top}\boldsymbol{W}^{-1}\boldsymbol{\xi}_{0}-\nu_{it}
\end{bmatrix},$$ and $\boldsymbol{W}_{t}(\boldsymbol{x})=\sum_{j=1}^{n}x_{j}\boldsymbol{W}_{tj}$, $\boldsymbol{R}_{t}(\boldsymbol{x})=\sum_{j=1}^{n}x_{j}\boldsymbol{r}_{tj}$, $H_{t}(\boldsymbol{x})=\sum_{j=1}^{n}x_{j}h_{tj}$.
\end{proposition}
\noindent $\textbf{Remark 6}$
\;\emph{We observe that by Propositions 2-4, $\mathit{\tilde{Z}}_{{C}_{2}}$ has a manifestly tractable representation in terms of Linear Matrix Inequalities (LMIs), which can be solved directly by using the powerful convex optimization solvers.}\\
\indent We prove that for joint binary DRCCP, i.e., $\mathit{S}\subseteq\left\{0,1\right\}^{n}$, the set $\mathit{S}\bigcap\mathit{Z}_{{C}_{2}}$ can be expressed as a mixed-integer convex reformulation when the constraint function is affine in both the decision variable and the uncertainty.
\begin{theorem}\label{11}
Suppose that $\mathit{S}\subseteq\left\{0,1\right\}^{n}$, $f_{t}\left(\boldsymbol{x},\boldsymbol{\xi}\right)=(\boldsymbol{A}^{t}\boldsymbol{x}+\boldsymbol{a}^{t})^{\top}\boldsymbol{\xi}+(\boldsymbol{b}^{t})^{\top}\boldsymbol{x}+h^{t}$ for any $t\in \left[T\right]$, and $\boldsymbol{\alpha}$ in $(14)$ can be upper bounded by a vector $\boldsymbol{M}$ for any $\boldsymbol{x}\in\mathit{S}$. Consider a convex set
\begin{equation}\notag
~~~~~~~~~~~\mathit{\hat{Z}}_{{C}_{2}}=
\left\{{\begin{array}{*{20}{c}}{\boldsymbol{x}\in\mathbb{R}^{n}:}&\begin{array}{l}
\lambda\delta+\frac{1}{N}\sum_{i=1}^{N}{s_{i}}\leq\epsilon,\\
G_{\hat{f}_{t}}\left(\boldsymbol{z}_{it},1,(\alpha_{t},\boldsymbol{y}^{t})\right)+1-\boldsymbol{z}_{it}^{\top}\zeta^{i}-s_{i}\leq 0 ,~\forall i\in\left[N\right],\forall t\in\left[T\right],\\
0\leq y_{r}^{t}\leq M_{t}x_{r}, \alpha_{t}-M_{t}(1-x_{r})\leq y_{r}^{t}\leq \alpha_{t},~\forall r\in\left[n\right],\forall t\in\left[T\right],\\
\left \|\boldsymbol{z}_{it}\right \|_{*}-\lambda\leq0 ,~\forall i\in\left[N\right],\forall t\in\left[T\right],\\
\lambda\geq 0,~\alpha_{t}\geq 0,~s_{i}\geq0,~\forall i\in\left[N\right],\forall t\in\left[T\right],
\end{array}
\end{array}} \right\}\begin{array}{*{20}{c}}
{~~~~~~~~\left( 28a \right)}\\
{~~~~~~~~\left( 28b \right)}\\
{~~~~~~~~\left( 28c \right)}\\
{~~~~~~~~\left( 28d \right)}\\
{~~~~~~~~\left( 28e \right)}
\end{array}
\end{equation}
where $\hat{f_{t}}\left((\alpha_{t},\boldsymbol{y}^{t}),\boldsymbol{\xi}\right)=(\boldsymbol{A}^{t}\boldsymbol{y}^{t}+\boldsymbol{a}^{t}\alpha_{t})^{\top}\boldsymbol{\xi}+(\boldsymbol{b}^{t})^{\top}\boldsymbol{y}^{t}+h^{t}\alpha_{t}$ for any $t\in \left[T\right]$, then
\begin{equation}\notag
\mathit{S}\bigcap{\mathit{\hat{Z}}}_{{C}_{2}}=\mathit{S}\bigcap \mathit{Z}_{{C}_{2}}\subseteq\mathit{S}\bigcap \mathit{Z}_{D}.
\end{equation}
\end{theorem}
\noindent{\it Proof}.~\;If each function $f_{t}\left(\boldsymbol{x},\boldsymbol{\xi}\right)=(\boldsymbol{A}^{t}\boldsymbol{x}+\boldsymbol{a}^{t})^{\top}\boldsymbol{\xi}+(\boldsymbol{b}^{t})^{\top}\boldsymbol{x}+h^{t}$, then for any $t\in \left[T\right]$ and $i\in \left[N\right]$, we have
\begin{equation}\notag
G_{f_{t}}\left(\boldsymbol{z}_{it},\alpha_{t},\boldsymbol{x}\right)=
\sup_{\boldsymbol{\xi}\in\mathrm{\Xi}}\left\{\boldsymbol{z}_{it}^{\top}\boldsymbol{\xi}
-\alpha_{t}\left[(\boldsymbol{A}^{t}\boldsymbol{x}+\boldsymbol{a}^{t})^{\top}\boldsymbol{\xi}+(\boldsymbol{b}^{t})^{\top}\boldsymbol{x}+h^{t}\right]\right\}.
\end{equation}
\indent Now, we define new variables $\boldsymbol{y}^{t}$ as $\boldsymbol{y}^{t}=\alpha_{t}\boldsymbol{x}$ for any $t\in \left[T\right]$.\\
\indent Since $\alpha_{t}\leq M_{t}$ for any $t\in \left[T\right]$, and hence by McCormick inequalities due to [25], we obtain
\begin{equation}\notag
0\leq y_{r}^{t}\leq M_{t}x_{r},~\alpha_{t}-M_{t}(1-x_{r})\leq y_{r}^{t}\leq \alpha_{t},
\end{equation}
which is exact for any $\boldsymbol{x}\subseteq\left\{0,1\right\}^{n}$. Thus, we have $\mathit{S}\bigcap{\mathit{\hat{Z}}}_{{C}_{2}}=\mathit{S}\bigcap \mathit{Z}_{{C}_{2}}$.
\qed
\noindent$\textbf{Remark 7}$
\;\emph{We note that by Theorem 3, to optimize over $\mathit{S}\bigcap \mathit{Z}_{{C}_{2}}$, we only need to optimize over $\mathit{S}\bigcap{\mathit{\hat{Z}}}_{{C}_{2}}$, which is
a mixed-integer convex set.}\\
\noindent$\textbf{Remark 8}$
\;\emph{In Theorem 3, we have assumed that one can find an upper bound $\boldsymbol{M}$, and indeed Corollary 2 also provides a sufficient condition for the existence of the vector $\boldsymbol{M}$.}\\
\indent We show how to find the proper upper bound $\boldsymbol{M}$ on the variable $\boldsymbol{\alpha}$ for any $\boldsymbol{x}\in\mathit{S}\bigcap\mathit{Z}_{{C}_{2}}$ in the following example.\\
\noindent\textbf{Example 1}\:~\:Let $\mathrm{\Xi}=\mathbb{R}^{m}$ in Theorem 3, then the set $\mathit{Z}_{{C}_{2}}$ can be rewritten as
\begin{equation}\notag
~~~~~~~~~~~~~\mathit{Z}_{{C}_{2}}=
\left\{{\begin{array}{*{20}{c}}{\boldsymbol{x}\in\mathbb{R}^{n}:}&\begin{array}{l}
\lambda\delta+\frac{1}{N}\sum_{i=1}^{N}{s_{i}}\leq\epsilon,\\
1-s_{i}\leq\alpha_{t}\left((\boldsymbol{A}^{t}\boldsymbol{x}+\boldsymbol{a}^{t})^{\top}\boldsymbol{\zeta}^{i}
+(\boldsymbol{b}^{t})^{\top}\boldsymbol{x}+h^{t}\right),~\forall i\in\left[N\right],\forall t\in\left[T\right],\\
\alpha_{t}\left \|\boldsymbol{A}^{t}\boldsymbol{x}+\boldsymbol{a}^{t}\right \|_{*}-\lambda\leq 0,~\forall t\in\left[T\right],\\
\lambda\!\geq 0,~\alpha_{t}\!\geq 0,~s_{i}\geq0,~\forall i\in\left[N\right],\forall t\in\left[T\right].
\end{array}
\end{array}} \right\}\begin{array}{*{20}{c}}
{~~~~~~~\left( 29a \right)}\\
{~~~~~~~\left( 29b \right)}\\
{~~~~~~~\left( 29c \right)}\\
{~~~~~~~\left( 29d \right)}
\end{array}
\end{equation}
\indent For any given $t\in\left[T\right]$, we add up all of constraints $(29\text{b})$ and obtain the inequality as below:
\begin{equation}\tag{30}
N-\sum_{i=1}^{N}{s_{i}}\leq\alpha_{t}
\left((\boldsymbol{A}^{t}\boldsymbol{x}+\boldsymbol{a}^{t})^{\top}\sum_{i=1}^{N}{\boldsymbol{\zeta}^{i}}+
N(\boldsymbol{b}^{t})^{\top}\boldsymbol{x}+Nh^{t}\right),
\end{equation}
\noindent furthermore, by constraint $(29\text{a})$, we note that
\begin{equation}\tag{31}
N+N(\lambda\delta-\epsilon)\leq\alpha_{t}
\left((\boldsymbol{A}^{t}\boldsymbol{x}+\boldsymbol{a}^{t})^{\top}
\sum_{i=1}^{N}{\boldsymbol{\zeta}^{i}}+N(\boldsymbol{b}^{t})^{\top}\boldsymbol{x}+Nh^{t}\right), \end{equation}
\noindent which implies that
\begin{equation}\tag{32}
N+N(\lambda\delta-\epsilon)\leq\alpha_{t}\left(\left\|(\boldsymbol{A}^{t})^{\top}
\sum_{i=1}^{N}{\boldsymbol{\zeta}^{i}}\right\|_{1}+(\boldsymbol{{a}}^{t})^{\top}\sum_{i=1}^{N}{\boldsymbol{\zeta}^{i}}+
N\left\|\boldsymbol{b}^{t}\right\|_{1}+Nh^{t}\right).
\end{equation}
\indent We observe that if $\mu\leq0$, then inequality $(32)$ can be rewritten as $N+N(\lambda\delta-\epsilon)\leq0$, which yields a contradiction to the fact that $\lambda\delta\geq0$, and hence we must have $\mu>0$, where $\mu=\left\|(\boldsymbol{A}^{t})^{\top}\sum_{i=1}^{N}{\boldsymbol{\zeta}^{i}}\right\|_{1}+
(\boldsymbol{{a}}^{t})^{\top}\sum_{i=1}^{N}{\boldsymbol{\zeta}^{i}}+N\left\|\boldsymbol{b}^{t}\right\|_{1}+Nh^{t}.$\\
\indent Since for any $i\in\left[N\right]$, $s_{i}\geq0$, then by constraints $(29\text{a})$ and $(29\text{c})$, we have $\alpha_{t}\leq\frac{\lambda}{\gamma}\leq\frac{\epsilon}{\delta\gamma}$,
where $\gamma=\min_{\boldsymbol{x}\in\left\{0,1\right\}^{n},\boldsymbol{A}^{t}
\boldsymbol{x}+\boldsymbol{a}^{t}\ne0}\left\|\boldsymbol{A}^{t}\boldsymbol{x}+\boldsymbol{a}^{t}\right\|_{*}$.\\
\indent Thus, we obtain $\frac{N+N(\lambda\delta-\epsilon)}{\mu}\leq\alpha_{t}\leq\frac{\epsilon}{\delta\gamma}$, and then for any $t\in\left[T\right]$, we let $M_{t}=\frac{\epsilon}{\delta\gamma}$.
\section{Numerical Results}
\indent In this section, we present two groups of numerical studies to demonstrate the computational effectiveness of the proposed formulations, one is to compare the out-of-sample performance of the approximation approach proposed in Theorem 2 with that of the sample average approximation (SAA) method in the context of a transportation decision problem, the other is to compare the approximation model proposed in Theorem 3 with exact Big-M model proposed in [13] so as to evaluate the optimality gap of the approximation model relative to the true optimality in the context of a binary knapsack problem.
\subsection{Distributionally robust chance-constrained transportation decision problem}
\indent We now consider a transportation decision problem, which is a classical problem that has been extensively studied in many literatures. The following example is an adaptation from [6]. Nonetheless, for completeness, we provide the following description of the problem. Given a set of facilities indexed by $\mathit{K}=\left\{k:k=1,\cdots,m\right\}$ and a set of customer locations indexed by $\mathit{J}=\left\{j:j=1,\cdots,n\right\}$. The total production quantity at each facility $k$ is $a_{k}$. Further assume that each facility has a normalized production capacity of $L_{high}$, i.e., $a_{k}\leq L_{high}$. Similarly, the total demand quantity from each customer location $j$ is $b_{j}$. Assume that the total production quantity equals to the total demand quantity, and they must be larger than a minimum demand quantity $L_{low}$, i.e., $\sum_{k=1}^{m}a_{k}=\sum_{j=1}^{n}b_{j}\geq L_{low}$. We aim to determine how to transport these products after all the ordered products are manufactured by the facilities in order to minimize the total cost. To formulate the problem, we define the uncertain parameter $\boldsymbol{\xi}=(\boldsymbol{\xi}_{1},\boldsymbol{\xi}_{2},\cdots,\boldsymbol{\xi}_{m})^{T}$ to be the vector of the unit transportation cost, where $\boldsymbol{\xi}_{k}\!=(\xi_{k1},\xi_{k2},\cdots,\xi_{kn}),
k=1,2,\cdots,m$, and $x_{kj}$ to be the volume of transport from point $k$ to $j$:
\begin{alignat}{2}\notag
f(\boldsymbol{\xi},\boldsymbol{a},\boldsymbol{b})=\min_{x_{kj}} \quad &\sum_{k=1}^{m}\sum_{j=1}^{n}\xi_{kj}x_{kj}\\ \notag
\mbox{s.t.}\quad
&\sum_{j=1}^{n}x_{kj}=a_{k},~\forall k\in\mathit{K},\\ \notag
&\sum_{k=1}^{m}x_{kj}=b_{j},~\forall j\in\mathit{J},\\ \notag
&x_{kj}\geq0,~\forall k\in\mathit{K},~\forall j\in\mathit{J}.\notag
\end{alignat}
\indent It should be noted that the unit transportation cost $\boldsymbol{\xi}$ supported on a rectangle of the form $\mathrm{\Xi}=\left\{\boldsymbol{\xi}\in\mathbb{R}^{mn}:0\leq\boldsymbol{\xi}\leq\boldsymbol{d}\right\}$ are only revealed at the second stage, then to ensure that the total transportation cost is low with high probability, the first-stage decision can be formulated as a DRCCP problem:
\begin{alignat}{2}\notag
\min_{\boldsymbol{a},\boldsymbol{b},z} \quad &z &\tag{33a}\\
\mbox{s.t.}\quad
&\inf_{\mathbb{P}\in\mathcal{P}}\mathbb{P}\left\{\boldsymbol{\xi}:f(\boldsymbol{\xi},\boldsymbol{a},\boldsymbol{b})\leq z\right\}\geq 1-\epsilon, &\tag{33b}\\
&\sum_{k=1}^{m}a_{k}=\sum_{j=1}^{n}b_{j}\geq L_{low}, &\tag{33c}\\
&0\leq a_{k}\leq L_{high},b_{j}\geq0,~\forall k\in\mathit{K},~\forall j\in\mathit{J}.&\tag{33d}
\end{alignat}
\indent Note that $z-f(\boldsymbol{\xi},\boldsymbol{a},\boldsymbol{b})$ is convex in the uncertainty $\boldsymbol{\xi}$ and concave in the decision variable $(\boldsymbol{a},\boldsymbol{b},z)$, which obviously satisfies assumption $(\mathbf{A1})$.\\
\indent Next, we review different approaches to construct $\mathcal{P}$ from $N$ sample data points $\left\{\boldsymbol{\zeta^{i}}\right\}_{i\in\left[N\right]}\subseteq\mathrm{\Xi}$ generated from the true distribution $\mathbb{P}_{\text{true}}$. Thus, for the proposed model $(33)$, we compare the out-of-sample performance of the distributionally robust approach based on 1-Wasserstein ball (denoted as DRW Model) with that of the classical sample average approximation (SAA) method (denoted as SAA Model).\\
\indent For DRW Model, we use the Wasserstein ambiguity set under assumption $(\mathbf{A3})$ with $\text{L}_{1}$-norm as distance metric. Then, problem $(33)$ can therefore be amenable to the worst-case CVaR approximation method discussed in Section 3. We thus obtain the approximation model by Theorem 2 as follows:
\begin{alignat}{2}\notag
\min_{\boldsymbol{a},\boldsymbol{b},\boldsymbol{x}^{i},\boldsymbol{y}^{i},\boldsymbol{v}_{i},q_{i},\alpha,z} \quad \quad &z &\tag{34a}\\
&\!\!\!\!\!\!\!\!\!\!\!\!\!\!\!\!\!\!\!\!\!\!\!\!\!\!\!\!\!\!\!\!\!\!\!\mbox{s.t.}\quad\;\delta\lvert v_{ir}\rvert+\frac{1}{N}\sum_{i=1}^{N}{q_{i}}\leq\epsilon\alpha,\,\forall i\in\left[N\right],\,\forall r\in\left[mn\right],&\tag{34b}\\
&\!\!\!\!\!\!\!\!\!\!\!\!\!\!\!\!\!\!\!\!\!(\boldsymbol{y}^{i})^{\top}\boldsymbol{d}-z+\alpha-
\boldsymbol{v}_{i}^{\top}
\boldsymbol{\zeta}^{i}-q_{i}\leq0,\,\forall i\in\left[N\right],&\tag{34c}\\
&\!\!\!\!\!\!\!\!\!\!\!\!\!\!\!\!\!\!\!\!\!
\boldsymbol{v}_{i}+\boldsymbol{x}^{i}\leq\boldsymbol{y}^{i},\,\forall i\in\left[N\right],&\tag{34d}\\
&\!\!\!\!\!\!\!\!\!\!\!\!\!\!\!\!\!\!\!\!\!\sum_{j=1}^{n}x_{kj}^{i}=a_{k},\,\sum_{k=1}^{m}x_{kj}^{i}=b_{j},
\,\forall k\in[m],\,\forall j\in[n],\,\forall i\in\left[N\right],&\tag{34e}\\
&\!\!\!\!\!\!\!\!\!\!\!\!\!\!\!\!\!\!\!\!\!\sum_{k=1}^{m}a_{k}=\sum_{j=1}^{n}b_{j}\geq L_{low},&\tag{34f}\\
&\!\!\!\!\!\!\!\!\!\!\!\!\!\!\!\!\!\!\!\!\!0\leq a_{k}\leq L_{high},b_{j}\geq0,x_{kj}^{i}\geq0,\boldsymbol{y}^{i}\geq0,q_{i}\geq0,\alpha>0,
\,\forall k\in[m],\,\forall j\in[n],\,\forall i\in\left[N\right].&\tag{34g}
\end{alignat}
\indent For SAA Model, we set $\mathcal{P}=\{\mathbb{P}_{\boldsymbol{\tilde{\zeta}}}\}$, which corresponds to a Wasserstein ball centered at the empirical distribution with the Wasserstein radius $\varepsilon=0$. Then, problem $(33)$ simply reduces to the corresponding SAA problem, which can be expressed as the following linear optimization problem
\begin{alignat}{2}\notag
\min_{\boldsymbol{a},\boldsymbol{b},z} \quad &z &\tag{35a}\\
\mbox{s.t.}\quad
&\frac{1}{N}\sum_{i=1}^{N}\mathbb{I}_{\left\{z-f(\boldsymbol{\zeta}^{i},\boldsymbol{a},\boldsymbol{b})\geq 0\right\}}\geq 1-\epsilon, &\tag{35b}\\
&\sum_{k=1}^{m}a_{k}=\sum_{j=1}^{n}b_{j}\geq L_{low}, &\tag{35c}\\
&0\leq a_{k}\leq L_{high},b_{j}\geq0,~\forall k\in[m],~\forall j\in[n],&\tag{35d}
\end{alignat}
where $\left\{\boldsymbol{\zeta^{i}}\right\}_{i\in\left[N\right]}\subseteq\mathrm{\Xi}$ are sample data points generated from the true distribution $\mathbb{P}_{\text{true}}$. Furthermore, according to [26], we can convert this into a large mixed-integer linear programming (MILP) as follows:
\begin{alignat}{2}\notag
\min_{\boldsymbol{a},\boldsymbol{b},\boldsymbol{x}^{i},\boldsymbol{s},z} \quad \quad &z &\tag{36a}\\
&\!\!\!\!\!\!\!\!\!\!\!\!\!\!\!\!\!\!\!\!\!\!\!\!\mbox{s.t.}\quad z-{\boldsymbol{\zeta}^{i}}^{\top}\boldsymbol{x}^{i}-M_{i}(s_{i}-1)\geq0,~\forall i\in\left[N\right],&\tag{36b}\\
&\!\!\!\!\!\!\!\!\!\!\!\frac{1}{N}\sum_{i=1}^{N}s_{i}\geq 1-\epsilon,&\tag{36c}\\
&\!\!\!\!\!\!\!\!\!\!\!\sum_{j=1}^{n}x_{kj}^{i}=a_{k},\,\sum_{k=1}^{m}x_{kj}^{i}=b_{j},\forall k\in[m],\,\forall j\in[n],\,\forall i\in\left[N\right],&\tag{36d}\\
&\!\!\!\!\!\!\!\!\!\!\!\sum_{k=1}^{m}a_{k}=\sum_{j=1}^{n}b_{j}\geq L_{low},&\tag{36e}\\
&\!\!\!\!\!\!\!\!\!\!\!\boldsymbol{s}\in\left\{0,1\right\}^{N}\!\!\!,0\leq a_{k}\leq L_{high},b_{j}\geq0,x_{kj}^{i}\geq0,\,\forall k\in[m],\,\forall j\in[n],\,\forall i\in\left[N\right],&\tag{36f}
\end{alignat}
where $M_{i}$ is a sufficiently large positive constant.\\
\indent In the subsequent tests, we set $m=4$, $n=6$, $L_{low}=L_{high}=2$, and $\boldsymbol{d}=8\boldsymbol{e}_{mn}$. We assume that the true distribution $\mathbb{P}_{\text{true}}$ of the unit transportation cost $\boldsymbol{\xi}$ is lognormal, which is designed in the following manner. That is, $\xi_{kj}=\exp(\widetilde{\xi}_{kj})$, where $\widetilde{\xi}_{kj},\,k\in[m],\,j\in[n]$, represent jointly normally distributed random variables
with the mean $\boldsymbol{\mu}$ drawn uniformly from $[0,1]^{mn}$ and the covariance matrix
$\boldsymbol{\mathrm{\Sigma}}=\text{Diag}(\boldsymbol{\sigma})\boldsymbol{C}\text{Diag}(\boldsymbol{\sigma})$, where $\boldsymbol{C}\in\mathbb{S}_{+}^{mn}$ is a random correlation matrix and $\boldsymbol{\sigma}=2\boldsymbol{e}_{mn}$ is the vector of standard deviations. Each problem instance is solved with CPLEX 12.10 using the YALMIP interface on a desktop with a 4.10 GHz processor and 32GB RAM.\\
\indent To assess the out-of-sample performance of different data-driven methods mentioned above, we conduct out-of-sample experiments to compare the reliability of the performance guarantees for these two models. Then we compute the optimal solution $(\boldsymbol{a}_{1}^{*},\boldsymbol{b}_{1}^{*},z_{1}^{*})$ of DRW Model by solving problem (34). Similarly, the optimal solution $(\boldsymbol{a}_{2}^{*},\boldsymbol{b}_{2}^{*},z_{2}^{*})$ of SAA Model is obtained by solving problem (36). Moreover, the out-of-sample performance is measured by the following chance constraint
\begin{alignat}{2}\notag
\mathbb{P}_{\text{true}}\left\{\boldsymbol{\xi}:z^{*}-f(\boldsymbol{\xi},\boldsymbol{a}^{*},\boldsymbol{b}^{*})\geq0 \right\}\geq 1-\epsilon,&\tag{37a}
\end{alignat}
which can be estimated at high accuracy using 20,000 test samples generated from $\mathbb{P}_{\text{true}}$ by solving another SAA problem
\begin{alignat}{2}\notag
\frac{1}{N}\sum_{i=1}^{N}\mathbb{I}_{\left\{z^{*}-
f(\boldsymbol{\zeta}^{i},\boldsymbol{a}^{*},\boldsymbol{b}^{*})\geq0\right\}}\geq 1-\epsilon. &\tag{37b}
\end{alignat}
\indent Thus, in the out-of-sample evaluation, for a given optimal solution $(\boldsymbol{a}^{*},\boldsymbol{b}^{*},z^{*})$ of DRW Model or SAA Model, checking whether it is feasible for chance constraint $(37\text{a})$ simplifies to calculating its reliability using $(37\text{b})$ and then comparing with the risk level $1-\epsilon$. In the following subsections 4.1.1 and 4.1.2, we study how the Wasserstein radius $\delta$ and the sample size $N$ affect the reliability of the optimal solutions for DRW Model and SAA Model.
\subsubsection{Impact of the Wasserstein radius $\delta$ on reliability}
\indent In this subsection, we compute the error bars between the 20\% and 80\% quantiles as well as the mean values of reliability for DRW Model under different radius $\delta\in\{0.01,0.04,0.07,0.10,0.13,0.16,0.19\}$ with the sample size $N=10$ and $160$, respectively, averaged across 20 independent random instances. Besides, we set the risk level to $\epsilon=0.10$.\\
\indent Figure 1 displays the reliability of the optimal solutions for DRW Model as a
function of $\delta$. We observe that the reliability of DRW Model increases as $\delta$ grows. This is because
the larger the radius $\delta$ is, the more distributions the Wasserstein ball $\mathcal{P}$ includes, and accordingly a larger uncertainty distribution family can capture the more uncertainty and can be a truer response for the all uncertainties. Moreover, it should be noted that from Figure 1a, as long as the proper radius is chosen, DRW Model is able to generate a feasible solution even if the size of the sample data set is very small. On the other hand, Figure 1b shows that when the sample data size is large enough, DRW Model can obtain a feasible solution even with the small radius.
\begin{figure}
\centering
\centerline{\includegraphics[height=13cm,width=18cm]{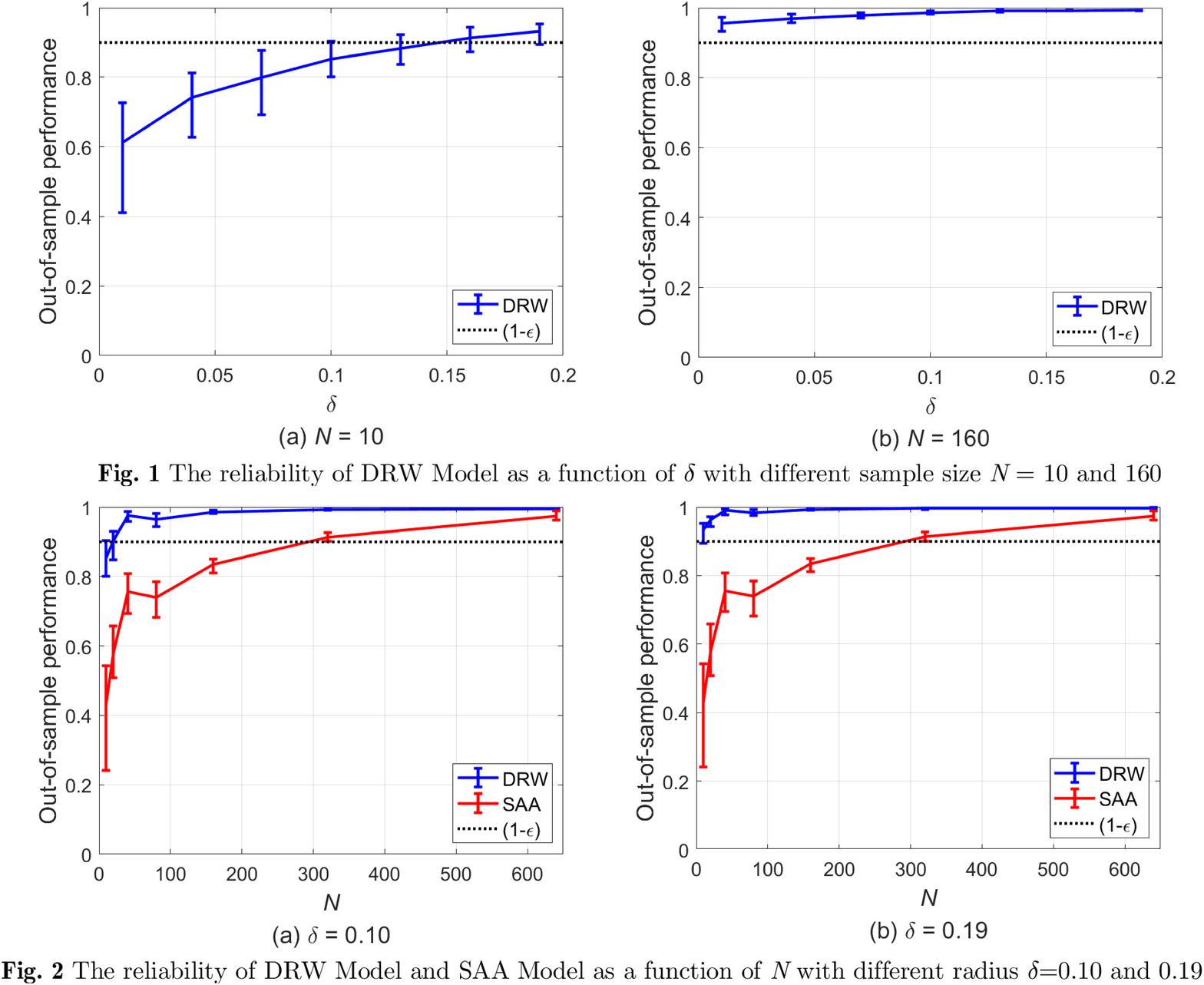}}
\label{1}
\end{figure}
\subsubsection{Impact of the sample size $N$ on reliability}
\indent In this subsection, we compute the error bars between the 20\% and 80\% quantiles as well as the mean values of reliability for DRW Model and SAA Model under different sample size $N\in\{10,20,40,80,160,320,640\}$ with the Wasserstein radius $\delta=0.10$ and $0.19$, respectively, averaged across 20 independent random instances. Besides, we set the risk level to $\epsilon=0.10$.\\
\indent Figure 2 depicts the reliability of the optimal solutions for DRW Model and SAA Model as
a function of $N$. We see that the reliability of both models tends to increase as the sample size becomes larger. Furthermore, Figure 2b shows that DRW Model can provide a high-quality reliable solution with a proper choice of radius $\delta$ even when the sample data points are very limited. However, SAA Model yields a poor reliability in situation where $N$ is small. In addition, the error bars visualize that DRW Model is significantly more stable than SAA Model. These results demonstrate that DRW Model is capable of returning reliable and stable solutions.
\subsection{Distributionally robust multidimensional knapsack problem}
\indent For the evaluation purpose, we study distributionally robust multidimensional knapsack problem (DRMKP) [27, 28] with binary decision variables. The following notations are adopted for binary DRMKP. We consider $T$ knapsacks and $n$ items, moreover, $c_{i}$ represents the value of item $i$ for any $i\in[n]$, $\boldsymbol{\xi}_{t}=\left(\xi_{t1},\cdots,\xi_{tn}\right)^{\top}$ represents the vector of random item weights supported on $\mathrm{\Xi}_{t}$ in knapsack $t$, and $b_{t}>0$ represents the capacity limit of knapsack $t$ for any $t\in[T]$. In addition, the decision variable $x_{i}\in\left\{0,1\right\}$ represents the proportion of $i$th item to be picked for any $i\in[n]$ and we let $\boldsymbol{x}\in\! \mathit{S}\!:=\left\{0,1\right\}^{n}$. Furthermore, we use the Wasserstein ambiguity sets under assumptions $(\mathbf{A2})$ and $(\mathbf{A3})$. With the notations above, binary DRMKP can be formulated as
\begin{alignat}{2}
\quad\quad\max_{\boldsymbol{x}} \quad &\boldsymbol{c}^{\top}\boldsymbol{x} &\tag{38a}\\
\mbox{s.t.}\quad
&\boldsymbol{x}\in\left\{0,1\right\}^{n}, &\tag{38b}\\
&\!\!\inf_{\mathbb{P}\in\mathcal{P}}\mathbb{P}\left\{\boldsymbol{\xi}\in\mathrm{\Xi}:
\boldsymbol{\xi}_{t}^{\top}\boldsymbol{x}\leq b_{t}, \forall t\in[T]\right\}\geq 1-\epsilon, &\tag{38c}
\end{alignat}
where constraint $(38\text{c})$ is to guarantee that the worst-case probability that the capacity of each knapsack should be satisfied is at least $1-\epsilon$.\\
\indent The following example demonstrates an application of Theorem 3. As observed in [13], the author derived an exact Big-M approach to solve it.\\
\noindent\textbf{Example 2}\:~\:Consider binary DRMKP $(38)$ with joint chance constraint $(T>1)$. Suppose the Wasserstein ambiguity set be defined as
\begin{alignat}{2}\notag
\mathcal{P}=\left\{\mathbb{P}:\mathbb{P}\left\{\boldsymbol{\xi}\in\mathrm{\Xi}\right\}=1,
\inf_{\mathbb{Q}}\left\{
\int_{\mathrm{\Xi}\times\mathrm{\Xi}}
\left\|\boldsymbol{\xi}-\boldsymbol{\tilde{\zeta}}\right\|_{2}\mathbb{Q}
\left(\mathit{d}\boldsymbol{\xi},\mathit{d}\boldsymbol{\tilde{\zeta}}\right)\right\}\leq\delta\right\},
\end{alignat}
where $\mathrm{\Xi}=\prod_{t\in[T]}\mathrm{\Xi}_{t}$ and $\mathrm{\Xi}_{t}=\mathbb{R}^{n}.$\\
\indent Then, by Theorem 3, binary DRMKP $(38)$ can be approximated by mixed-integer second-order cone programming (MISOCP) as follows:
\begin{alignat}{2}
\max\quad&\boldsymbol{c}^{\top}\boldsymbol{x} &\tag{39a}\\
\mbox{s.t.}\quad
&\boldsymbol{x}\in\left\{0,1\right\}^{n}, &\tag{39b}\\
&\lambda\delta+\frac{1}{N}\sum_{i=1}^{N}{s_{i}}\leq\epsilon, &\tag{39c}\\
&\!\!-\alpha_{t}b_{t}+(\boldsymbol{y}^{t})^{\top}\boldsymbol{\zeta}^{i}_{t}+1-s_{i}\leq 0,~\forall i\in\left[N\right],\forall t\in\left[T\right], &\tag{39d}\\
&\!0\leq y_{r}^{t}\leq \frac{\epsilon}{\delta}x_{r}, \alpha_{t}-\frac{\epsilon}{\delta}(1-x_{r})\leq y_{r}^{t}\leq \alpha_{t},~\forall r\in\left[n\right],\forall t\in\left[T\right], &\tag{39e}\\
&\!\!\left \|\boldsymbol{y}^{t}\right \|_{2}-\lambda\leq0,~\forall t\in\left[T\right], &\tag{39f}\\
&\!\lambda\geq 0,\alpha_{t}\geq 0,s_{i}\geq0,~\forall i\in\left[N\right],\forall t\in\left[T\right].&\tag{39g}
\end{alignat}
\indent It should be noted that according to Example 1, $\alpha_{t}$ can be upper bounded by $M_{t}=\frac{\epsilon}{\delta}$ for all $t\in\left[T\right]$ in problem $(39)$.\\
\indent By Example 2, we present one group of numerical experiments, i.e., to compare our approximation method proposed in Theorem 3 with exact Big-M method proposed in [13].\\
\indent The problem instances are created as follows. We generate 20 random instances with $n=20$ and $T=10$. For each instance, we generate $N\in\left\{100,1000\right\}$ empirical samples $\left\{\boldsymbol{\zeta}^{i}\right\}_{i\in\left[N\right]}\subseteq\mathbb{R}_{+}^{n\times T}$ from a uniform distribution over a box $\left[1,10\right]^{n\times T}$. For each $i\in\left[n\right]$, we independently generate $c_{i}$ from the uniform distribution on the interval $\left[1,10\right]$. In addition, we set $b_{t}=50$ for each $t\in\left[T\right]$. We test these 20 random instances with the risk level $\epsilon\in\left\{0.05,0.10\right\}$ and the Wasserstein radius $\delta\in\left\{0.01, 0.02\right\}$. All the instances are executed on a desktop with a 4.10 GHz processor and 32GB RAM, while CPLEX 12.10 is used with their default settings. We set the time limit of solving each instance to be 3600 seconds.\\
\indent The numerical results with sample size $N = 100$ are displayed in Table 1. We use Big-M Model to denote exact reformulation proposed in [13]. Besides, we use CVaR Model to denote inner approximation proposed in Theorem 3. We use ``Opt.Val'' to denote the optimal value $v^{*}$, ``Value'' to denote the best objective value output from the approximation model $(39)$, and ``Time'' to denote the total running time in seconds. Additionally, since we can solve Big-M Model to the optimality, we use ``GAP'' to denote the optimality gap of the approximation model $(39)$, which is computed as
$$\text{GAP}=\frac{\left|\text{Value}-\text{Opt.Val}\right|}{\text{Opt.Val}}.$$
\indent The numerical results with sample size $N = 1000$ are displayed in Table 2. Similarly, we use Big-M Model to denote exact reformulation proposed in [13]. Besides, we use CVaR Model to denote inner approximation proposed in Theorem 3. We use ``LB'' to denote the best lower bounds found by the models that we're going to test, and ``Time'' to denote the total running time in seconds. To evaluate the effectiveness of CVaR Model, we use ``Improvement'' to denote the percentage of differences between the lower bound of the approximation model $(39)$ and the lower bound of Big-M Model, which is computed as
$$\text{Improvement}=\frac{\text{LB of Approximation Model}-\text{LB of Exact Model}}{\text{LB of Exact Model}},$$
where Approximation Model here represents the approximation model $(39)$. Similarly, Exact Model here represents Big-M Model.\\
\indent From Table 1, we note that CVaR Model can be solved to the optimality within 15 seconds, while Big-M Model often takes longer time to solve. In terms of approximation accuracy, CVaR Model nearly finds the true optimal solutions in most instances. These results demonstrate that CVaR Model is capable of finding near-optimal solutions.\\
\indent From Table 2, we observe that the total running time of CVaR Model significantly outperforms that of Big-M Model, i.e., the majority of the instances of CVaR Model can be solved within 20 minutes, while Big-M Model cannot be solved within the time limit. In terms of approximation accuracy, we see that CVaR Model can find at least the same feasible solutions as Big-M Model. Additionally, in some instances, CVaR Model can provide slightly better lower bounds than Big-M Model. These results demonstrate the effectiveness of CVaR Model, which scales for large number of sample data points.\\
\indent As can be seen from Table 1 and Table 2, the main reasons for the different numerical performances between Big-M Model and CVaR Model are as follows: $(i)$ CVaR Model only involves $\mathcal{O}(n)$ binary variables, while Big-M Model involves $\mathcal{O}(N+n)$ binary variables. $(ii)$ Big-M Model requires more auxiliary variables than CVaR Model.
\begin{table}[H]
\caption{Numerical results of Big-M Model proposed in [13] and CVaR Model proposed in Theorem 3 for binary DRMKP when sample size N = 100.}
\begin{center}
\begin{spacing}{1.10}
\scalebox{0.75}{
\begin{tabular}{ccccc|cc|cccc}
\hline
\multicolumn{3}{c|}{\multirow{2}*{$T$~~~~~~~~~~$n$~~~~~~~~~~~$\epsilon$~~~~~~~~~~~~$\delta$~~~~~~~~~~~~~Instances }}&~~~~~Big-M\;\;\;\;\;Model&&~~~~~~$\rm{CVaR}$\;\;\;\;\;Model\\
\cline{4-8}
\multicolumn{3}{c|}{~} &\!\!\!\!\!\!\!\!\!\!\!\!\!\!\!\!\!\!\!Opt.Val&\!\!\!\!\!\!\!\!\!\!\!\!\!\!\!\!\!\!\!\!Time&\!\!\!\!\!\!\!\!\!\!\!\!\!\!\!\!\!\!\!\!\!\!\!\!\!\!\!\!Value&\!\!\!\!\!\!\!\!\!\!\!\!\!\!\!\!\!\!\!\!\!\!\!\!\!\!\!\!\!GAP&~~\!\!\!~\!Time\\
\hline
&&~~~~~~~~~~~~~~1&\!\!\!\!\!\!\!\!\!\!\!\!\!\!\!\!\!\!50.49&\!\!\!\!\!\!\!\!\!\!\!\!\!\!\!\!\!\!\!\!4.82&\!\!\!\!\!\!\!\!\!\!\!\!\!\!\!\!\!\!\!\!\!\!\!\!\!\!\!50.49&\!\!\!\!\!\!\!\!\!\!\!\!\!\!\!\!\!\!\!\!\!\!\!\!\!\!\!0.00\%&~12.75\\
&&~~~~~~~~~~~~~~2&\!\!\!\!\!\!\!\!\!\!\!\!\!\!\!\!\!\!50.32&\!\!\!\!\!\!\!\!\!\!\!\!\!\!\!\!\!\!\!\!2.32&\!\!\!\!\!\!\!\!\!\!\!\!\!\!\!\!\!\!\!\!\!\!\!\!\!\!\!50.32&\!\!\!\!\!\!\!\!\!\!\!\!\!\!\!\!\!\!\!\!\!\!\!\!\!\!\!0.00\%&~6.29\\
&&~~~~~~~~~~~~~~3&\!\!\!\!\!\!\!\!\!\!\!\!\!\!\!\!\!\!49.51&\!\!\!\!\!\!\!\!\!\!\!\!\!\!\!\!\!\!\!\!7.67&\!\!\!\!\!\!\!\!\!\!\!\!\!\!\!\!\!\!\!\!\!\!\!\!\!\!\!49.51&\!\!\!\!\!\!\!\!\!\!\!\!\!\!\!\!\!\!\!\!\!\!\!\!\!\!\!0.00\%&~8.02\\
&&~~~~~~~~~~~~~~4&\!\!\!\!\!\!\!\!\!\!\!\!\!\!\!\!\!\!50.10&\!\!\!\!\!\!\!\!\!\!\!\!\!\!\!\!\!\!\!\!9.67&\!\!\!\!\!\!\!\!\!\!\!\!\!\!\!\!\!\!\!\!\!\!\!\!\!\!\!49.90&\!\!\!\!\!\!\!\!\!\!\!\!\!\!\!\!\!\!\!\!\!\!\!\!\!\!\!0.40\%&~13.34\\
&&~~~~~~~~~~~~~~5&\!\!\!\!\!\!\!\!\!\!\!\!\!\!\!\!\!\!49.55&\!\!\!\!\!\!\!\!\!\!\!\!\!\!\!\!\!\!\!\!2.18&\!\!\!\!\!\!\!\!\!\!\!\!\!\!\!\!\!\!\!\!\!\!\!\!\!\!\!48.25&\!\!\!\!\!\!\!\!\!\!\!\!\!\!\!\!\!\!\!\!\!\!\!\!\!\!\!2.63\%&~9.52\\
&&~~~~~~~~~~~~~~6&\!\!\!\!\!\!\!\!\!\!\!\!\!\!\!\!\!\!44.00&\!\!\!\!\!\!\!\!\!\!\!\!\!\!\!\!\!\!\!\!3.07&\!\!\!\!\!\!\!\!\!\!\!\!\!\!\!\!\!\!\!\!\!\!\!\!\!\!\!43.96&\!\!\!\!\!\!\!\!\!\!\!\!\!\!\!\!\!\!\!\!\!\!\!\!\!\!\!0.07\%&~7.66\\
\hline
\end{tabular}}
\end{spacing}
\end{center}
\end{table}

\begin{table}[H]
\caption{Numerical results of Big-M Model proposed in [13] and CVaR Model proposed in Theorem 3 for binary DRMKP when sample size N = 1000.}
\begin{center}
\begin{spacing}{1.10}
\scalebox{0.75}{
\begin{tabular}{ccccc|cc|cccc}
\hline
\multicolumn{3}{c|}{\multirow{2}*{$T$~~~~~$n$~~~~~~$\epsilon$~~~~~~~$\delta$~~~~~~~~Instances }}&~~~~~Big-M\;\;\;\;\;Model&&~~~~~~~~~~~~~~$\rm{CVaR}$\;\;\;\;\;Model\\
\cline{4-8}
\multicolumn{3}{c|}{~} &\!\!\!\!\!\!\!\!\!\!\!\!\!\!\!\!\!\!\!\!\!\!\!\!\!LB&\!\!\!\!\!\!\!\!\!\!\!\!\!\!\!\!\!\!\!\!\!\!\!Time&\!\!\!\!\!\!\!\!\!\!\!\!\!\!\!\!\!\!\!\!\!\!\!\!\!\!\!\!\!\!\!\!\!\!\!\!\!\!\!\!\!\!\!~~~~LB&\!\!\!\!\!\!\!\!\!\!\!\!\!\!\!\!\!\!\!\!\!\!\!\!\!\!\!\!\!\!\!\!\!\!\!\!\!\!\!\!\!Improvement&\!Time\\
\hline
&&\!\!\!\!\!\!\!\!\!\!1&\!\!\!\!\!\!\!\!\!\!\!\!\!\!\!\!\!\!\!\!\!\!\!\!\!40.44&\!\!\!\!\!\!\!\!\!\!\!\!\!\!\!\!\!\!\!\!\!\!3600&\!\!\!\!\!\!\!\!\!\!\!\!\!\!\!\!\!\!\!\!\!\!\!\!\!\!\!\!\!\!\!\!\!\!\!40.44&\!\!\!\!\!\!\!\!\!\!\!\!\!\!\!\!\!\!\!\!\!\!\!\!\!\!\!\!\!\!\!\!\!\!\!\!\!\!\!\!\!0.00\%&1171.59\\
&&\!\!\!\!\!\!\!\!\!\!2&\!\!\!\!\!\!\!\!\!\!\!\!\!\!\!\!\!\!\!\!\!\!\!\!\!40.18&\!\!\!\!\!\!\!\!\!\!\!\!\!\!\!\!\!\!\!\!\!\!3600&\!\!\!\!\!\!\!\!\!\!\!\!\!\!\!\!\!\!\!\!\!\!\!\!\!\!\!\!\!\!\!\!\!\!\!41.14&\!\!\!\!\!\!\!\!\!\!\!\!\!\!\!\!\!\!\!\!\!\!\!\!\!\!\!\!\!\!\!\!\!\!\!\!\!\!\!\!\!2.40\%&1032.97\\
&&\!\!\!\!\!\!\!\!\!\!3&\!\!\!\!\!\!\!\!\!\!\!\!\!\!\!\!\!\!\!\!\!\!\!\!\!43.70&\!\!\!\!\!\!\!\!\!\!\!\!\!\!\!\!\!\!\!\!\!\!3600&\!\!\!\!\!\!\!\!\!\!\!\!\!\!\!\!\!\!\!\!\!\!\!\!\!\!\!\!\!\!\!\!\!\!\!43.70&\!\!\!\!\!\!\!\!\!\!\!\!\!\!\!\!\!\!\!\!\!\!\!\!\!\!\!\!\!\!\!\!\!\!\!\!\!\!\!\!\!0.00\%&476.74\\
&&\!\!\!\!\!\!\!\!\!\!4&\!\!\!\!\!\!\!\!\!\!\!\!\!\!\!\!\!\!\!\!\!\!\!\!\!43.65&\!\!\!\!\!\!\!\!\!\!\!\!\!\!\!\!\!\!\!\!\!\!3600&\!\!\!\!\!\!\!\!\!\!\!\!\!\!\!\!\!\!\!\!\!\!\!\!\!\!\!\!\!\!\!\!\!\!\!43.65&\!\!\!\!\!\!\!\!\!\!\!\!\!\!\!\!\!\!\!\!\!\!\!\!\!\!\!\!\!\!\!\!\!\!\!\!\!\!\!\!\!0.00\%&1375.80\\
\hline
\hline
\end{tabular}}
\end{spacing}
\end{center}
\end{table}
\section{Conclusion}
\indent In this paper, we studied distributionally robust joint chance-constrained programming with convex uncertain constraints under Wasserstein ambiguity set. We proposed an equivalent reformulation of the set $\mathit{{Z}_{D}}$ and then developed the worst-case CVaR approximation of the set $\mathit{{Z}_{D}}$ via a system of biconvex constraints, which is naturally hard to solve. Furthermore, a convex relaxation of the proposed approximation can be derived by constructing new decision variables which allows us to eliminate biconvex terms. Once the decision variables are binary and the uncertain
constraints are affine, this proposed biconvex approximation is equivalent to a tractable mixed-integer convex programming. Numerical results demonstrated that the proposed models can be solved efficiently. In our study, we assume the uncertain mapping $f_{t}\left(\boldsymbol{x},\boldsymbol{\xi}\right)$ is convex in $\boldsymbol{\xi}$ and concave in $\boldsymbol{x}$. A future direction is to consider DRCCP problems with a broader family of uncertain mappings, for instance, when each constraint function is quasi-convex in the uncertainty and is concave in the decision variable.\\
\noindent \textbf{Data Availability Statement}\\
\indent Some or all data, models, or code generated or used during the study are available from the first author of the paper by request.\\
\noindent \textbf{Reference}\\
\noindent1. Calafiore, G.C., ElGhaoui, L.: On distributionally robust chance-constrained linear programs. J. Optim. Theory Appl. \textbf{130}(1), 1-22 (2006)\\
\noindent2. Hanasusanto, G.A., Roitch, V., Kuhn, D., Wiesemann, W.: A distributionally robust perspective on uncertainty quantification and chance constrained programming. Math. Program. \textbf{151}, 35-62 (2015)\\
\noindent3. Hanasusanto, G.A., Roitch, V., Kuhn, D., Wiesemann, W.: Ambiguous joint chance constraints under mean and dispersion information. Oper. Res. \textbf{65}(3), 751-767 (2017)\\
\noindent4. Jiang, R., Guan, Y.: Data-driven chance constrained stochastic program. Math. Program. \textbf{158}, 291-327 (2016)\\
\noindent5. Xie, W., Ahmed, S.: On deterministic reformulations of distributionally robust joint chance constrained optimization problems. SIAM J. Optim. \textbf{28}(2), 1151-1182 (2018)\\
\noindent6. Yang, W., Xu, H.: Distributionally robust chance constraints for non-linear uncertainties. Math. Program. \textbf{155}, 231-265 (2016)\\
\noindent7. Chen, W., Sim, M., Sun, J., Teo, C.P.: From CVaR to uncertainty set: Implications in joint chance constrained optimization. Oper. Res. \textbf{58}(2), 470-485 (2010)\\
\noindent8. Delage, E. and Ye, Y.: Distributionally robust optimization under moment uncertainty with application to data-driven problems. Oper. Res. \textbf{58}(3), 595-612 (2010)\\
\noindent9. Zymler, S., Kuhn, D., Rustem, B.: Distributionally robust joint chance constraints with second-order moment information. Math. Program. \textbf{137}, 167-198 (2013)\\
\noindent10. Wiesemann, W., Kuhn, D., Sim, M.: Distributionally robust convex optimization. Oper. Res. \textbf{62}(6), 1358-1376 (2014)\\
\noindent11. Xie, W., Ahmed, S., Jiang, R.: Optimized Bonferroni approximations of distributionally robust joint chance constraints. Math. Program. (2019) https://doi.org/10.1007/s10107-019-01442-8\\
\noindent12. Chen, Z., Kuhn, D., Wiesemann, W.: Data-driven chance constrained programs over Wasserstein balls. Available on Optimization Online (2018) http://www.optimization-online.org/DB\_HTML/2018/06/66\\
71.html\\
\noindent13. Xie, W.: On distributionally robust chance constrained programs with Wasserstein distance. Math. Program. \textbf{186}, 115-155 (2021)\\
\noindent14. Ji, R., Lejeune, M.A.: Data-driven distributionally robust chance-constrained optimization with Wasserstein metric. J Glob Optim. (2020) https://doi.org/10.1007/s10898-020-00966-0\\
\noindent15. Hota, A.R., Cherukuri, A., Lygeros, J.: Data-driven chance constrained optimization under Wasserstein ambiguity sets. In 2019 American Control Conference (ACC), pp. 1501-1506, IEEE, (2019) https://doi.org/10.23919/ACC.2019.8814677\\
\noindent16. Ho-Nguyen, N., K\i l\i n\c{c}-Karzan, F., K\"{u}\c{c}\"{u}kyavuz, S., Lee, D.: Distributionally robust chance-constrain-\\
ed programs with right-hand side uncertainty under Wasserstein ambiguity. Math. Program. (2021) https://doi.org/10.1007/s10107-020-01605-y\\
\noindent17. Ho-Nguyen, N., K\i l\i n\c{c}-Karzan, F., K\"{u}\c{c}\"{u}kyavuz, S., Lee, D.: Strong formulations for distributionally robust chance-constrained programs with left-hand side uncertainty under Wasserstein ambiguity. (2020) arXiv:2007.06750v2\\
\noindent18. Esfahani, P.M., Kuhn, D.: Data-driven distributionally robust optimization using the Wasserstein metric: Performance guarantees and tractable reformulations. Math. Program. \textbf{171}(1-2), 115-166 (2018)\\
\noindent19. Gao, R., Kleywegt, A.J.: Distributionally robust stochastic optimization with Wasserstein distance. Accepted by Math. Oper. Res. (2016) https://doi.org/10.48550/arXiv.1604.02199\\
\noindent20. Mehrotra, S., Papp, D.: A cutting surface algorithm for semi-infinite convex programming with an application to moment robust optimization. SIAM J. Optim. \textbf{24}(4), 1670-1697 (2014)\\
\noindent21. Luo, F., Mehrotra, S.: Decomposition algorithm for distributionally robust optimization using Wasserstein metric with an application to a class of regression models. European J. Oper. Res. \textbf{278}(1), 20-35 (2019)\\
\noindent22. Rockafellar, R.T., Uryasev, S.: Optimization of conditional value-at-risk. J. Risk \textbf{2}, 21-42 (2000)\\
\noindent23. Shapiro, A., Kleywegt, A.: Minimax analysis of stochastic problems. Optim. Methods Softw. \textbf{17}(3), 523-542 (2002)\\
\noindent24. Pichler, A., Xu, H.: Quantitative stability analysis for minimax distributionally robust risk optimization. Math. Program. (2018) https://doi.org/10.1007/s10107-018-1347-4\\
\noindent25. Mccormick, G.P.: Computability of global solutions to factorable nonconvex programs: Part I-convex underestimating problems. Math. Program. \textbf{10}(1), 147-175 (1976)\\
\noindent26. Ruszczy$\acute{\text{n}}$ski, A.: Probabilistic programming with discrete distributions and precedence constrained knapsack polyhedra. Math. Program. \textbf{93}, 195-215 (2002)\\
\noindent27. Cheng, J., Delage, E., Lisser, A.: Distributionally robust stochastic knapsack problem. SIAM J. Optim. \textbf{24}(3), 1485-1506 (2014)\\
\noindent28. Song, Y., Luedtke, J.R., K\"{u}\c{c}\"{u}kyavuz, S.: Chance-constrained binary packing problems. INFORMS J. Comput. \textbf{26}(4), 735-747 (2014)\\
\noindent \textbf{Appendix A: Proofs}\\
\noindent \textbf{Proof of Proposition 2.} Note that $(24\text{b})$ can be rewritten as
\begin{alignat}{2}\notag
q_{it}\geq &G_{f_{t}}\left(\boldsymbol{v}_{it},1,\boldsymbol{x}\right)+
\alpha_{t}-\boldsymbol{v}_{it}^{\top}\boldsymbol{\zeta^{i}}\\ \notag
&=\sup_{\boldsymbol{\xi}\in\mathrm{\Xi}}\left[\boldsymbol{v}_{it}^{\top}\boldsymbol{\xi}
-\left(\boldsymbol{\xi}^{\top}\boldsymbol{x}+\langle\boldsymbol{A}^{t},\boldsymbol{\xi}\boldsymbol{\xi}^{\top}\rangle+(\boldsymbol{b}^{t})^{\top}\boldsymbol{x}+h^{t}\right)\right]+\alpha_{t}-\boldsymbol{v}_{it}^{\top}\boldsymbol{\zeta^{i}}\\ \notag
&=\sup_{\boldsymbol{\xi}\in\mathrm{\Xi}}\left[(\boldsymbol{v}_{it}-\boldsymbol{x})^{\top}\boldsymbol{\xi}-\langle\boldsymbol{A}^{t},\boldsymbol{\xi}\boldsymbol{\xi}^{\top}\rangle\right]-\left[(\boldsymbol{b}^{t})^{\top}\boldsymbol{x}+h^{t}\right]+\alpha_{t}-\boldsymbol{v}_{it}^{\top}\boldsymbol{\zeta^{i}}. \tag{40}
\end{alignat}
\indent We now discuss the reformulation of the inner subproblem $\sup_{\boldsymbol{\xi}\in\mathrm{\Xi}}$ in (40). For this problem, there exists $\boldsymbol{\xi}$ such that $\boldsymbol{a}_{k}^{\top}\boldsymbol{\xi}< d_{k}$, $k=1,2,\cdots,l$ due to $d_{k}>0$, which implies Slater's condition is satisfied, and hence the strong duality holds.\\
\indent Thus, we have
$$\sup_{\boldsymbol{\xi}\in\mathrm{\Xi}}\left[(\boldsymbol{v}_{it}-\boldsymbol{x})^{\top}\boldsymbol{\xi}-
\langle\boldsymbol{A}^{t},\boldsymbol{\xi}\boldsymbol{\xi}^{\top}\rangle\right]=
\inf_{\boldsymbol{\nu^{it}}\geq0}\sup_{\boldsymbol{\xi}\in\mathbb{R}^{m}}
\left[(\boldsymbol{v}_{it}-\boldsymbol{x})^{\top}\boldsymbol{\xi}-
\langle\boldsymbol{A}^{t},\boldsymbol{\xi}\boldsymbol{\xi}^{\top}\rangle
-\sum_{k=1}^{l}\nu_{k}^{it}(\boldsymbol{a}_{k}^{\top}\boldsymbol{\xi}-d_{k})\right],$$
which is equivalent to
\begin{alignat}{2}\notag
\min_{\boldsymbol{\nu^{it}}\geq0,u_{it}\in\mathbb{R}} \quad &u_{it}\\ \notag
\mbox{s.t.}\quad
&(\boldsymbol{v}_{it}-\boldsymbol{x})^{\top}\boldsymbol{\xi}-\langle\boldsymbol{A}^{t},\boldsymbol{\xi}\boldsymbol{\xi}^{\top}\rangle
-\sum_{k=1}^{l}\nu_{k}^{it}(\boldsymbol{a}_{k}^{\top}\boldsymbol{\xi}-d_{k})\leq u_{it},~\forall \boldsymbol{\xi}\in\mathbb{R}^{m},\notag
\end{alignat}
which can be further written as
\begin{alignat}{2}\notag
\min_{\boldsymbol{\nu^{it}}\geq0,u_{it}\in\mathbb{R}} \quad &u_{it}\\ \notag
\mbox{s.t.}\quad
&\begin{bmatrix}
\boldsymbol{A}^{t}&-\frac{1}{2}(\boldsymbol{v}_{it}-\boldsymbol{x}-\sum_{k=1}^{l}\nu_{k}^{it}\boldsymbol{a}_{k})\\
-\frac{1}{2}(\boldsymbol{v}_{it}-\boldsymbol{x}-\sum_{k=1}^{l}
\nu_{k}\boldsymbol{a}_{k})^{\top}&u_{it}-\sum_{k=1}^{l}\nu_{k}^{it}d_{k}
\end{bmatrix}\succeq0,\notag
\end{alignat}
and thus, $(40)$ is equivalent to
\begin{alignat}{2}\notag
q_{it}\geq\min_{\boldsymbol{\nu^{it}}\geq0,u_{it}\in\mathbb{R}} &u_{it}-\left[(\boldsymbol{b}^{t})^{\top}\boldsymbol{x}+h^{t}\right]+\alpha_{t}-\boldsymbol{v}_{it}^{\top}\boldsymbol{\zeta^{i}}\\
\notag
\mbox{s.t.}\quad
&\begin{bmatrix}
\boldsymbol{A}^{t}&-\frac{1}{2}(\boldsymbol{v}_{it}-\boldsymbol{x}-\sum_{k=1}^{l}\nu_{k}^{it}\boldsymbol{a}_{k})\\
-\frac{1}{2}(\boldsymbol{v}_{it}-\boldsymbol{x}-\sum_{k=1}^{l}\nu_{k}\boldsymbol{a}_{k})^{\top}&u_{it}-\sum_{k=1}^{l}\nu_{k}^{it}d_{k}
\end{bmatrix}\succeq0.\notag
\end{alignat}
Then the claim follows.
\qed
\noindent \textbf{Proof of Proposition 3.} Note that similar to the proof of Proposition 2, $(24\text{b})$ can be rewritten as
\begin{alignat}{2}\notag
q_{it}\geq &G_{f_{t}}\left(\boldsymbol{v}_{it},1,\boldsymbol{x}\right)+\alpha_{t}
-\boldsymbol{v}_{it}^{\top}\boldsymbol{\zeta^{i}}\\ \notag
&=\sup_{\boldsymbol{\xi}\in\mathrm{\Xi}}\left[(\boldsymbol{v}_{it}-\boldsymbol{x})^{\top}\boldsymbol{\xi}-\langle\boldsymbol{A}^{t},\boldsymbol{\xi}\boldsymbol{\xi}^{\top}\rangle\right]-\left[(\boldsymbol{b}^{t})^{\top}\boldsymbol{x}+h^{t}\right]+\alpha_{t}-\boldsymbol{v}_{it}^{\top}\boldsymbol{\zeta^{i}}.
\tag{41}
\end{alignat}
\indent We now discuss the reformulation of the inner subproblem $\sup_{\boldsymbol{\xi}\in\mathrm{\Xi}}$ in (41). For this problem, there exists $\boldsymbol{\xi}=\boldsymbol{\xi}_{0}$ such that $({\boldsymbol{\xi}-\boldsymbol{\xi}_{0}})^{\top}\boldsymbol{W}^{-1}({\boldsymbol{\xi}-\boldsymbol{\xi}_{0}})<1$, which implies Slater's condition is satisfied, and hence the strong duality holds.\\
\indent Thus, we have
\begin{alignat}{2}\notag
&\sup_{\boldsymbol{\xi}\in\mathrm{\Xi}}\left[(\boldsymbol{v}_{it}-
\boldsymbol{x})^{\top}\boldsymbol{\xi}-
\langle\boldsymbol{A}^{t},\boldsymbol{\xi}\boldsymbol{\xi}^{\top}\rangle\right]\\
&=\inf_{\nu_{it}\geq0}\sup_{\boldsymbol{\xi}\in\mathbb{R}^{m}}\left[(\boldsymbol{v}_{it}-
\boldsymbol{x})^{\top}\boldsymbol{\xi}-\langle\boldsymbol{A}^{t},\boldsymbol{\xi}
\boldsymbol{\xi}^{\top}\rangle
-\nu_{it}(({\boldsymbol{\xi}-\boldsymbol{\xi}_{0}})^{\top}\boldsymbol{W}^{-1}
({\boldsymbol{\xi}-\boldsymbol{\xi}_{0}})-1)\right],\notag
\end{alignat}
which is equivalent to
\begin{alignat}{2}\notag
\min_{\nu_{it}\geq0,u_{it}\in\mathbb{R}} \quad &u_{it}\\ \notag
\mbox{s.t.}\quad
&(\boldsymbol{v}_{it}-\boldsymbol{x})^{\top}\boldsymbol{\xi}-
\langle\boldsymbol{A}^{t},\boldsymbol{\xi}\boldsymbol{\xi}^{\top}\rangle-
\nu_{it}(({\boldsymbol{\xi}-\boldsymbol{\xi}_{0}})^{\top}\boldsymbol{W}^{-1}({\boldsymbol{\xi}-\boldsymbol{\xi}_{0}})-1)\leq u_{it},~\forall \boldsymbol{\xi}\in\mathbb{R}^{m},\notag
\end{alignat}
which can be further written as
\begin{alignat}{2}\notag
\min_{\nu_{it}\geq0,u_{it}\in\mathbb{R}} \quad &u_{it}\\ \notag
\mbox{s.t.}\quad
&\begin{bmatrix}
\boldsymbol{A}^{t}+\nu_{it}\boldsymbol{W}^{-1}&-\frac{1}{2}(2\nu_{it}\boldsymbol{W}^{-1}\boldsymbol{\xi}_{0}+
\boldsymbol{v}_{it}-\boldsymbol{x})\\
-\frac{1}{2}(2\nu_{it}\boldsymbol{W}^{-1}\boldsymbol{\xi}_{0}+\boldsymbol{v}_{it}-
\boldsymbol{x})^{\top}&u_{it}+\nu_{it}\boldsymbol{\xi}_{0}^{\top}\boldsymbol{W}^{-1}\boldsymbol{\xi}_{0}-\nu_{it}
\end{bmatrix}\succeq0,\notag
\end{alignat}
and thus, $(41)$ is equivalent to
\begin{alignat}{2}\notag
q_{it}\geq\min_{\nu_{it}\geq0,u_{it}\in\mathbb{R}} &u_{it}-\left[(\boldsymbol{b}^{t})^{\top}\boldsymbol{x}+h^{t}\right]+
\alpha_{t}-\boldsymbol{v}_{it}^{\top}\boldsymbol{\zeta^{i}}\\
\notag
\mbox{s.t.}\quad
&\begin{bmatrix}
\boldsymbol{A}^{t}+\nu_{it}\boldsymbol{W}^{-1}&
-\frac{1}{2}(2\nu_{it}\boldsymbol{W}^{-1}\boldsymbol{\xi}_{0}+\boldsymbol{v}_{it}-\boldsymbol{x})\\
-\frac{1}{2}(2\nu_{it}\boldsymbol{W}^{-1}\boldsymbol{\xi}_{0}+\boldsymbol{v}_{it}-
\boldsymbol{x})^{\top}&u_{it}+\nu_{it}\boldsymbol{\xi}_{0}^{\top}\boldsymbol{W}^{-1}\boldsymbol{\xi}_{0}-\nu_{it}
\end{bmatrix}\succeq0.\notag
\end{alignat}
Then the claim follows.
\qed
\noindent \textbf{Proof of Proposition 4.}
Note that $(24\text{b})$ can be rewritten as
\begin{alignat}{2}\notag
q_{it}\geq &G_{f_{t}}\left(\boldsymbol{v}_{it},1,\boldsymbol{x}\right)+
\alpha_{t}-\boldsymbol{v}_{it}^{\top}\boldsymbol{\zeta^{i}}\\ \notag
&=\sup_{\boldsymbol{\xi}\in\mathrm{\Xi}}\left[\boldsymbol{v}_{it}^{\top}\boldsymbol{\xi}
-\boldsymbol{w_{t}}(\boldsymbol{\xi})^{\top}\boldsymbol{x}\right]+
\alpha_{t}-\boldsymbol{v}_{it}^{\top}\boldsymbol{\zeta^{i}}\\ \notag
&=\sup_{\boldsymbol{\xi}\in\mathrm{\Xi}}\left[\boldsymbol{v}_{it}^{\top}\boldsymbol{\xi}
-\left(\boldsymbol{\xi}^{\top}\boldsymbol{W}_{t}(\boldsymbol{x})\boldsymbol{\xi}+
\boldsymbol{R}_{t}(\boldsymbol{x})^{\top}\boldsymbol{\xi}+H_{t}(\boldsymbol{x})\right)\right]
+\alpha_{t}-\boldsymbol{v}_{it}^{\top}\boldsymbol{\zeta^{i}}\\ \notag
&=\sup_{\boldsymbol{\xi}\in\mathrm{\Xi}}\left[(\boldsymbol{v}_{it}-
\boldsymbol{R}_{t}(\boldsymbol{x}))^{\top}\boldsymbol{\xi}-\boldsymbol{\xi}^{\top}
\boldsymbol{W}_{t}(\boldsymbol{x})\boldsymbol{\xi}\right]-H_{t}(\boldsymbol{x})+\alpha_{t}-
\boldsymbol{v}_{it}^{\top}\boldsymbol{\zeta^{i}}.\tag{42}
\end{alignat}
\indent We now discuss the reformulation of the inner subproblem $\sup_{\boldsymbol{\xi}\in\mathrm{\Xi}}$ in (42). Similar to the proof of Proposition 3, for this problem, there exists $\boldsymbol{\xi}=\boldsymbol{\xi}_{0}$ such that $({\boldsymbol{\xi}-\boldsymbol{\xi}_{0}})^{\top}\boldsymbol{W}^{-1}({\boldsymbol{\xi}-\boldsymbol{\xi}_{0}})<1$, which implies Slater's condition is satisfied, and hence the strong duality holds.\\
\indent Thus, we have
\begin{alignat}{2}\notag
&\sup_{\boldsymbol{\xi}\in\mathrm{\Xi}}\left[(\boldsymbol{v}_{it}-\boldsymbol{R}_{t}(\boldsymbol{x}))^{\top}\boldsymbol{\xi}-\boldsymbol{\xi}^{\top}\boldsymbol{W}_{t}(\boldsymbol{x})\boldsymbol{\xi}\right]\\\notag
&=\inf_{\nu_{it}\geq0}\sup_{\boldsymbol{\xi}\in\mathbb{R}^{m}}\left[(\boldsymbol{v}_{it}-\boldsymbol{R}_{t}(\boldsymbol{x}))^{\top}\boldsymbol{\xi}-\boldsymbol{\xi}^{\top}\boldsymbol{W}_{t}(\boldsymbol{x})\boldsymbol{\xi}-\nu_{it}(({\boldsymbol{\xi}-\boldsymbol{\xi}_{0}})^{\top}\boldsymbol{W}^{-1}({\boldsymbol{\xi}-\boldsymbol{\xi}_{0}})-1)\right],
\end{alignat}
which is equivalent to
\begin{alignat}{2}\notag
\min_{\nu_{it}\geq0,u_{it}\in\mathbb{R}} \quad &u_{it}\\ \notag
\mbox{s.t.}\quad
&(\boldsymbol{v}_{it}-\boldsymbol{R}_{t}(\boldsymbol{x}))^{\top}\boldsymbol{\xi}-
\boldsymbol{\xi}^{\top}\boldsymbol{W}_{t}(\boldsymbol{x})\boldsymbol{\xi}-\nu_{it}(({\boldsymbol{\xi}-\boldsymbol{\xi}_{0}})^{\top}\boldsymbol{W}^{-1}({\boldsymbol{\xi}-\boldsymbol{\xi}_{0}})-1)\leq u_{it},~\forall \boldsymbol{\xi}\in\mathbb{R}^{m},\notag
\end{alignat}
which can be further written as
\begin{alignat}{2}\notag
\min_{\nu_{it}\geq0,u_{it}\in\mathbb{R}} \quad &u_{it}\\ \notag
\mbox{s.t.}\quad
&\begin{bmatrix}
\boldsymbol{W}_{t}(\boldsymbol{x})+\nu_{it}\boldsymbol{W}^{-1}&-\frac{1}{2}(2\nu_{it}\boldsymbol{W}^{-1}\boldsymbol{\xi}_{0}+\boldsymbol{v}_{it}-\boldsymbol{R}_{t}(\boldsymbol{x}))\\
-\frac{1}{2}(2\nu_{it}\boldsymbol{W}^{-1}\boldsymbol{\xi}_{0}+\boldsymbol{v}_{it}-\boldsymbol{R}_{t}(\boldsymbol{x}))^{\top}&u_{it}+\nu_{it}\boldsymbol{\xi}_{0}^{\top}\boldsymbol{W}^{-1}\boldsymbol{\xi}_{0}-\nu_{it}
\end{bmatrix}\succeq0,\notag
\end{alignat}
and thus, $(42)$ is equivalent to
\begin{alignat}{2}\notag
q_{it}\geq\min_{\nu_{it}\geq0,u_{it}\in\mathbb{R}} &u_{it}-H_{t}(\boldsymbol{x})+\alpha_{t}-\boldsymbol{v}_{it}^{\top}\boldsymbol{\zeta^{i}}\notag\\
\mbox{s.t.}\quad
&\begin{bmatrix}
\boldsymbol{W}_{t}(\boldsymbol{x})+\nu_{it}\boldsymbol{W}^{-1}&-\frac{1}{2}(2\nu_{it}\boldsymbol{W}^{-1}\boldsymbol{\xi}_{0}+\boldsymbol{v}_{it}-\boldsymbol{R}_{t}(\boldsymbol{x}))\\
-\frac{1}{2}(2\nu_{it}\boldsymbol{W}^{-1}\boldsymbol{\xi}_{0}+\boldsymbol{v}_{it}-
\boldsymbol{R}_{t}(\boldsymbol{x}))^{\top}&u_{it}+\nu_{it}\boldsymbol{\xi}_{0}^{\top}\boldsymbol{W}^{-1}
\boldsymbol{\xi}_{0}-\nu_{it}
\end{bmatrix}\succeq0.\notag
\end{alignat}
Then the claim follows.
\qed
\end{document}